\input amstex
\input epsf
\magnification=\magstep1 
\baselineskip=13pt
\documentstyle{amsppt}
\vsize=8.7truein
\CenteredTagsOnSplits \NoRunningHeads

\def\AA{\Cal A}
\def\var{\bold{var\thinspace}}
\def\cov{\bold{cov\thinspace}}
\def\EE{\bold{E\thinspace }}
\def\Pr{\bold{Pr\thinspace }}

\def\vl{\operatorname{vol}}
\def\rk{\operatorname{rank}}
\def\conv{\operatorname{conv}}

\title Maximum entropy Gaussian approximation for the number of integer points and 
volumes of polytopes   \endtitle
\author Alexander Barvinok and John Hartigan \endauthor
\address Department of Mathematics, University of Michigan, Ann Arbor,
MI 48109-1043, USA \endaddress
\email barvinok$\@$umich.edu \endemail
\address Department of Statistics, Yale University, New Haven, CT 06520-8290 \endaddress
\email john.hartigan$\@$yale.edu \endemail 
\thanks The research of the first author was partially supported by NSF Grant DMS 0856640 and a
United States - Israel BSF grant 2006377. \endthanks
 \abstract  We describe a maximum entropy approach for computing
 volumes and counting integer points in polyhedra.  
 To estimate the number of points from a particular set 
 $X \subset {\Bbb R}^n$ in a polyhedron $P \subset {\Bbb R}^n$, by solving a certain entropy maximization problem,
  we construct a probability distribution on the set $X$  such that a) the probability mass function is constant on the set $P \cap X$ and b) the expectation 
 of the distribution lies in $P$. This allows us to apply Central Limit Theorem 
 type arguments to deduce computationally efficient  approximations for the number 
 of integer points, volumes, and the number of 0-1 vectors in the polytope. 
 As an application, we obtain asymptotic formulas for volumes of multi-index 
 transportation polytopes and for the number of multi-way contingency tables.
 \endabstract
\date July 2009 
\enddate
\endtopmatter
\document

\head 1. Introduction  \endhead 

In this paper, we address the problems of computing the volume and counting the number of integer 
points in a given polytope. These problems have a long history, see for example, surveys 
\cite{GK94}, \cite{DL05} and \cite{Ve05},
and, generally speaking, are computationally hard. We describe a maximum entropy approach which, in a 
number of non-trivial cases, allows one to obtain good quality approximations by solving certain 
specially constructed convex optimization problems on 
polytopes. Those optimization problems can be solved quite efficiently, in theory and in practice,
by interior point methods, see \cite{NN94}. 

The essence of our approach is as follows: given a discrete set $S \subset {\Bbb R}^n$ of interest, such as 
the set ${\Bbb Z}^n_+$ of all non-negative integer points or the set $\{0, 1\}^n$ of all 0-1 points, 
and an affine subspace $A \subset {\Bbb R}^n$ we want to compute or estimate the number 
$|S \cap A|$ of points in $A$. For that, we construct a probability measure $\mu$ on 
$S$ with the property that the probability mass function is constant on the set $A \cap S$ and 
the expectation of $\mu$ lies in $A$. These two properties allow us to apply Local Central Limit 
Theorem type arguments to estimate $|S \cap A|$. The measure $\mu$ turns out to be the measure 
of the largest entropy on $S$ with the expectation in $A$, so that constructing $\mu$ 
reduces to solving a convex optimization problem. We also consider a continuous version 
of the problem, where $S$ is the non-negative orthant ${\Bbb R}^n_+$ and 
our goal is to estimate the volume of the set $S \cap A$. 

Our approach is similar in spirit 
to that of E.T. Jaynes \cite{Ja57} (see also \cite{Go63}), who, motivated by problems of statistical mechanics, 
formulated a general principle of estimating the average value of a functional $g$ 
with respect to an unknown probability distribution  on a discrete set $S$ of states provided the average values of some other functionals $f_1, \ldots, f_r$ on $S$ are given. He suggested to estimate $g$ by its expectation with respect to the maximum entropy probability distribution on $S$ such that the expectations of $f_i$ have prescribed values. Our situation fits this general framework when,
for example, $S$ is the set
${\Bbb Z}^n_+$ of non-negative integer vectors, $f_i$ are the equations defining an affine 
subspace $A$,  functional $g$ is some quantity of interest, while the unknown probability distribution on $S$ is the counting measure on $S \cap A$ 
(in interesting cases, the set $S \cap A$ is complicated enough so that we may justifiably 
think of the counting measure on $S \cap A$ as of an unknown measure).

\subhead (1.1) Definitions and notation \endsubhead 
In what follows, ${\Bbb R}^n$ is Euclidean space with the standard integer lattice 
${\Bbb Z}^n \subset {\Bbb R}^n$. A polyhedron $P \subset {\Bbb R}^n$ is defined as the 
set of solutions $x=\left(\xi_1, \ldots, \xi_n \right)$ to a vector equation
$$\xi_1 a_1 + \ldots + \xi_n a_n =b, \tag1.1.1$$
where $a_1, \ldots, a_n; b \in {\Bbb R}^d$ are $d$-dimensional vectors for $d<n$, 
and inequalities 
$$\xi_1, \ldots, \xi_n \ \geq \ 0. \tag1.1.2$$
We assume that vectors $a_1, \ldots, a_n$ span ${\Bbb R}^d$, in which case the affine 
subspace defined by (1.1.1) has dimension $n-d$. We also assume that $P$ has a non-empty 
interior, that is, contains a point $x=\left(\xi_1, \ldots, \xi_n \right)$, where inequalities (1.1.2) 
are strict. One of our goals is to compute the $(n-d)$-dimensional volume $\vl P$ of $P$ with respect to the 
Lebesgue measure in the affine subspace (1.1.1) induced from ${\Bbb R}^n$. More generally, 
our approach allows us to estimate the {\it exponential integral} 
$$\int_P e^{\ell(x)} \ dx,$$
where $\ell: {\Bbb R}^n \longrightarrow {\Bbb R}$ is a linear function. We note that the integral may be 
well defined even if $P$ is unbounded. Often, we use a shorthand $Ax=b, x \geq 0$ for (1.1.1)--(1.1.2),
where $A=\left[ a_1, \ldots, a_n \right]$ is the matrix with the columns $a_1, \ldots, a_n$ and $x$ is 
thought of as a column vector $x=\left[ \xi_1, \ldots, \xi_n \right]^T$.

We are also interested in the number $|P \cap {\Bbb Z}^n|$ of integer points in $P$. In this case, we assume 
that vectors $a_1, \ldots, a_n$ and $b$ are integer, that is, $a_1, \ldots, a_n; b \in {\Bbb Z}^d$. 
The number $|P \cap {\Bbb Z}^n|$ as a function of vector $b$ in (1.1.1) is known as 
the {\it vector partition function} associated with vectors $a_1, \ldots, a_n$, see 
for example, \cite{BV97}. More generally, 
our approach allows us to estimate the {\it exponential sum}
$$\sum_{m \in P \cap {\Bbb Z}^n} e^{\ell(m)},$$
where $\ell: {\Bbb R}^n \longrightarrow {\Bbb R}$ is a linear function. Again, the sum may converge 
even if polyhedron $P$ is unbounded.

Finally, we consider a version of the integer point counting problem where we are interested 
in 0-1 vectors only. Namely, let $\{0, 1\}^n$ be the set (Boolean cube) of all vectors in ${\Bbb R}^n$
with the coordinates 0 and 1. We estimate $|P \cap \{0, 1\}^n|$ and, more generally, the sum
$$\sum_{m \in P \cap \{0, 1\}^n} e^{\ell(m)}.$$

\subhead (1.2) The maximum entropy approach \endsubhead Let us consider 
the integer counting problem first. One of the most straightforward approaches to computing 
$|P \cap {\Bbb Z}^n|$ approximately is via the Monte Carlo method. As in Section 1.1, we think of $P$ as 
defined by a system $Ax=b, x \geq 0$. 
 We place $P$ in a sufficiently 
large axis-parallel integer box $B$ in the non-negative orthant ${\Bbb R}^n_+$ of 
${\Bbb R}^n$, sample integer points from $B$ independently at random and 
count what proportion of points lands in $P$. It is well understood that the method is very inefficient
if $P$ occupies a small fraction of $B$, in which case the sampled points will not land in $P$ unless we 
use great many samples. Let $X$ be a random vector distributed 
uniformly on the set of integer points in box $B$. One can try to circumvent sampling entirely by considering 
the random vector $Y=AX$ and interpreting the number of integer points in $P$ in terms of the 
probability mass function of $Y$ at $b$. One can hope then, in the spirit of the Central Limit Theorem,
 that since the coordinates of $Y$ are 
linear combinations of independent coordinates $x_1, \ldots, x_n$ of $X$, the 
distribution of $Y$ is somewhat close to the Gaussian and hence the probability mass function 
of $Y$ at $b$ can be approximated by the Gaussian density. The problem with this approach is 
that, generally speaking, the expectation $\EE Y$ will be very far from the target vector $b$, so 
one tries to apply the Local  Central Limit Theorem on the tail of the distribution, which is 
precisely where it is not applicable.

We propose a simple remedy to this naive Monte Carlo approach. Namely, by solving a convex 
optimization problem on $P$, we construct a multivariate 
geometric random variable $X$ such that 
\bigskip
(1.2.1) The probability mass function of $X$ is constant on the set $P \cap {\Bbb Z}^n$ of integer points in $P$;
\medskip
(1.2.2) We have $\EE X \in P$, or, equivalently,  $\EE Y =b$ for $Y=AX$.
\bigskip
Condition (1.2.1) allows us to express the number $|P \cap {\Bbb Z}^n|$ of integer points in $P$ 
in terms of the probability mass function of $Y$, while condition (1.2.2) allows us to prove 
the Local Central Limit Theorem for $Y$ in a variety of situations. We have 
$X=\left(x_1, \ldots, x_n \right)$ where $x_j$ are independent geometric random variables 
with expectations 
$\zeta_j$ such that $z=\left(\zeta_1, \ldots, \zeta_n\right)$ is the unique point maximizing the value 
of the strictly concave function, the entropy of $X$,
$$g(x)=\sum_{j=1}^n \Bigl( \left(\xi_j +1 \right) \ln \left( \xi_j +1 \right) - \xi_j \ln \xi_j \Bigr)$$
on $P$, see Theorem 3.1 for the precise statement.

Similarly, to estimate the number of 0-1 vectors in $P$, we construct a multivariate  Bernoulli 
random variable $X$, such that (1.2.2) holds while (1.2.1) is replaced by 
\bigskip
(1.2.3) The probability mass function of $X$ is constant on the set $P \cap \{0, 1\}^n$ of 0-1 
vectors in $P$.
\bigskip

In this case, $X=\left(x_1, \ldots, x_n \right)$, where $x_j$ are independent Bernoulli random 
variables with expectations $\zeta_j$ such that
$z=\left(\zeta_1, \ldots, \zeta_n \right)$ is the unique point maximizing the value of the
strictly concave function, the entropy of $X$,
$$h(x)=\sum_{j=1}^n \left(\xi_j \ln {1 \over \xi_j} + (1-\xi_j) \ln {1 \over 1-\xi_j} \right)$$
on the truncated polytope 
$$P \cap \Bigl\{ 0 \leq \xi_j \leq 1: \quad \text{for} \quad j=1, \ldots, n \Bigr\},$$
see Theorem 3.3 for the precise statement.

Finally, to approximate the volume of $P$, we construct a multivariate exponential random 
variable $X$ such that (1.2.2) holds and (1.2.1) is naturally replaced 
by
\bigskip
(1.2.4) The density of $X$ is constant on $P$.
\bigskip
Condition (1.2.4) allows us to express the volume of $P$ in terms of the density of $Y=AX$ at $Y=b$, while (1.2.2) allows us to establish a Local Central Limit Theorem for $Y$ in a number of cases. 
In this case, each coordinate $x_j$ is sampled independently from the exponential 
distribution with expectation $\zeta_j$ such that $z=\left(\zeta_1, \ldots, \zeta_n \right)$ is the 
unique point maximizing the value of the strictly concave function, the entropy of $X$,
$$f(x)=n+\sum_{j=1}^n \ln \xi_j $$
on $P$, see Theorem 3.6 for the precise statement.  In optimization, the point $z$ is known as the {\it analytic center} of $P$ and it played a 
central role in the development of interior point methods, see \cite{Re88}.

These three examples (counting integer points, counting 0-1 vectors, and computing volumes) 
are important particular cases of a general approach to counting through the solution to an 
entropy maximization problem (cf. Theorem 3.5) with the subsequent asymptotic analysis 
of multivariate integrals needed to establish the Local Central Limit Theorem type results.

\head 2. Main results \endhead 

\subhead (2.1) Gaussian approximation for volume \endsubhead 
Let $P \subset {\Bbb R}^n$ be a polytope,
defined by a system $Ax=b, x \geq 0$, where $A$ is an $d \times n$ matrix with the columns 
$a_1, \ldots, a_n$. We assume that
$\rk A=d<n$.
We find the point $z=\left(\zeta_1, \ldots, \zeta_n \right)$ maximizing 
$$f(x)=n+\sum_{j=1}^n \ln \xi_j, \qquad x=\left(\xi_1, \ldots, \xi_n \right)$$ 
on $P$. Let $B$ be the $d \times n$ matrix with the columns $\zeta_1 a_1, \ldots, \zeta_n a_n$.
We approximate the volume of $P$ by the Gaussian formula
$$\vl P \approx {1 \over (2 \pi)^{d/2}} \left({\det A A^T \over \det B B^T} \right)^{1/2} e^{f(z)} \tag2.1.1$$
We consider the standard scalar product $\langle \cdot, \cdot \rangle$ and the 
corresponding Euclidean norm $\| \cdot \|$ in ${\Bbb R}^d$.

We prove the following main result. 
\proclaim{(2.2) Theorem} Let us consider a quadratic form 
$q: {\Bbb R}^d \longrightarrow {\Bbb R}$ defined by 
$$q(t)={1 \over 2} \sum_{j=1}^n \zeta_j^2 \langle a_j, t \rangle^2.$$
Suppose that for some $\lambda>0$ we have
$$q(t) \geq \lambda \|t\|^2 \quad \text{for all} \quad t \in {\Bbb R}^d$$ 
and that for some $\theta >0$ we have 
$$\zeta_j \|a_j\| \leq \theta \quad \text{for} \quad j=1, \ldots, n.$$
Then there exists an absolute constant $\gamma$ such that the following holds:
\break let $0 < \epsilon  \leq 1/2$ be a number and suppose that 
$$\lambda \ \geq \ \gamma \theta^2 \epsilon^{-2} \left( d+ \ln {1 \over \epsilon}\right)^2 
\ln  \left({n \over \epsilon}\right).$$
Then the number 
$${1 \over (2 \pi)^{d/2}} \left({\det A A^T \over \det B B^T} \right)^{1/2} e^{f(z)}$$ 
approximates $\vl P$ within relative error $\epsilon$.
\endproclaim

Let us consider the columns $a_1, \ldots, a_n$ of $A$ as vectors from Euclidean space 
${\Bbb R}^d$ endowed with the standard scalar product $\langle \cdot, \cdot \rangle$.
The quadratic form $q$ defines the moment of inertia of the set of vectors 
$\left\{\zeta_1 a_1, \ldots, \zeta_n a_n\right\}$,
see, for example, \cite{Ba97}. By requiring that the smallest eigenvalue of $q$ is sufficiently large 
compared to the lengths of  the vectors $\zeta_j a_j$, we require that the set is sufficiently 
``round''. For a sufficiently generic (random) set of $n$ vectors, we will have $q(t)$ 
roughly proportional to $\|t\|^2$ and hence $\lambda$ will be of the order of 
$n d^{-1} \max_{j=1, \ldots, n} \zeta_j^2 \|a_j\|^2$. 

We prove Theorem 2.2 in Section 6.

In Section 4, we apply Theorem 2.2 to approximate the volume of a multi-index transportation polytope,
see, for example, \cite{Y+84}, that is, the polytope $P$ of $\nu$-dimensional 
$k_1 \times \ldots \times k_{\nu}$ arrays of 
non-negative numbers $\left(\xi_{j_1 \ldots j_{\nu}} \right)$ with $1 \leq j_i \leq k_i$ for 
$i=1, \ldots, \nu$ with prescribed sums along coordinate hyperplanes $j_i=j$.
 We show that Theorem 2.2 implies that asymptotically the volume of $P$ is given by a Gaussian formula (2.1.1)
 as long as $\nu \geq 5$. We suspect that the Gaussian approximation holds as long as 
 $\nu \geq 3$, but the proof would require some additional considerations beyond those of Theorem 2.2.
 In particular, for $\nu \geq 5$ we obtain the asymptotic formula for the volume of the 
 polytope of polystochastic tensors, see \cite{Gr92}. 

For $\nu=2$ polytope $P$ is the usual transportation polytope. 
Interestingly, its volume is {\it not} given by the Gaussian formula, cf. \cite{CM07b}. 

 In \cite{Ba09}, a much cruder asymptotic formula  of $\vl P$ in terms of 
$e^{f(z)}$ was proved under much weaker assumptions.

\subhead (2.3) Gaussian approximation for the number of integer points \endsubhead
For a polytope $P$, defined by a system $Ax=b, x \geq 0$, we find the point 
$z=\left(\zeta_1, \ldots, \zeta_n\right)$ maximizing 
$$g(x)=\sum_{j=1}^n \Bigl( \left(\xi_j +1 \right) \ln \left( \xi_j +1 \right) - \xi_j \ln \xi_j \Bigr),
\qquad x=\left(\xi_1, \ldots, \xi_n \right)$$
on $P$. Assuming that $a_1, \ldots, a_n \in {\Bbb Z}^d$ are the columns of $A$, we define $B$ as the 
$d \times n$ matrix whose $j$-th column is $\left(\zeta_j^2 +\zeta_j\right)^{1/2} a_j$ for $j=1, \ldots, n$.

 We assume that $A$ is an integer $d \times n$ matrix of rank $d < n$. Let 
 $\Lambda=A\left({\Bbb Z}^n \right)$ be image of the standard lattice, $\Lambda \subset {\Bbb Z}^d$.
 We approximate the number of integer points in $P$ by the Gaussian formula
 $$ \left| P \cap {\Bbb Z}^n \right| \approx {e^{g(z)}\det \Lambda  \over (2 \pi)^{d/2} (\det BB^T)^{1/2}}.
 \tag2.3.1$$
 In this paper, we consider the simplest case of $\Lambda={\Bbb Z}^d$,
  which is equivalent to the 
 greatest common divisor of the $d \times d$ minors of $A$ being equal to 1.
 
 Together with the Euclidean norm $\| \cdot \|$ in ${\Bbb R}^d$, we consider the $\ell^1$ and 
$\ell^{\infty}$ norms:
$$\|t\|_1=\sum_{i=1}^d |\tau_i| \quad \text{and} \quad \|t\|_{\infty} =
\max_{i=1, \ldots, d} |\tau_i| \quad \text{where} \quad t=\left(\tau_1, \ldots, \tau_d \right).$$ 
Clearly, we have 
$$\|t\|_1 \ \geq \ \|t\| \ \geq \ \|t\|_{\infty} \quad \text{for all} \quad t \in {\Bbb R}^d.$$
Compared to the case of volume estimates (Sections 2.1--2.2), we acquire an additive error which is 
governed by the arithmetic of the problem. 

Let $e_1, \ldots, e_d$ be the standard 
basis of ${\Bbb Z}^d$.
 We prove the 
 following main result.
 
 \proclaim{(2.4) Theorem} Let us consider a quadratic form $q: {\Bbb R}^d \longrightarrow {\Bbb R}$
defined by 
$$q(t)={1 \over 2} \sum_{j=1}^n \left(\zeta_j +\zeta_j^2 \right) \langle a_j, t \rangle^2.$$
For $i=1, \ldots, d$ let us choose a non-empty finite set $Y_i \subset {\Bbb Z}^n$ such that 
$Ay=e_i$ for all $y \in Y_i$ and let us define a quadratic form 
$\psi_i: {\Bbb R}^n \longrightarrow {\Bbb R}$ by 
$$\psi_i(x)={1 \over |Y_i|} \sum_{y \in Y_i} \langle y, x \rangle^2.$$
Suppose that for some $\lambda \geq 0$ we have 
$$q(t) \geq \lambda \|t\|^2 \quad \text{for all} \quad t \in {\Bbb R}^d,$$
that for some $\rho > 0$ we have 
$$\psi_i(x) \leq \rho \|x\|^2 \quad \text{for all} \quad x \in {\Bbb R}^n \quad \text{and} \quad i=1, \ldots, d,$$
 that for some $\theta \geq 1$ we have 
$$ \| a_j \|_1  \leq \theta \sqrt{\zeta_j \over \left(1+\zeta_j\right)^3} \quad \text{for} \quad j=1, \ldots, n$$
and that
$$\zeta_j (1+\zeta_j) \geq \alpha \quad \text{for} \quad j=1, \ldots, n$$
and some $\alpha \geq 0$.

Then, for some absolute constant $\gamma>0$ and for any $0 \leq \epsilon \leq 1/2$,
as long as 
$$\lambda \geq \gamma  \epsilon^{-2} \theta^2 \left(d + \ln {1\over \epsilon} \right)^2
\ln \left({n \over \epsilon}\right),$$
we have 
$$\left| P \cap {\Bbb Z}^n \right| =
e^{g(z)} \left( {\kappa \over (2 \pi)^{d/2} \left(\det BB^T \right)^{1/2}} +\Delta \right),$$
where 
$$1-\epsilon \leq \kappa \leq 1+\epsilon$$
and
$$|\Delta| \leq  \left(1 +{2 \over 5} \alpha \pi^2 \right)^{-m} \quad \text{for} \quad
m=\left\lfloor {1 \over 16 \pi^2 \rho  \theta^2} \right\rfloor.$$
\endproclaim

While the condition on the smallest eigenvalue of quadratic form $q$ is very similar to that of 
Theorem 2.2 and is linked to the metric properties of $P$, the appearance of quadratic forms 
$\psi_i$ is explained by the arithmetic features of $P$. Let us choose $1 \leq i \leq d$ and let 
us consider the affine subspace $\AA_i$ of the points $x \in {\Bbb R}^n$ such that $Ax=e_i$.
Let $\Lambda_i =\AA_i \cap {\Bbb Z}^n$ be the point lattice in $\AA_i$. We would like to choose 
a set $Y_i \subset \Lambda_i$ in such a way that the maximum eigenvalue $\rho_i$ of the form 
$\psi_i$, which defines the moment of inertia of $Y_i$, see \cite{Ba97}, becomes as small as 
possible, $\rho_i \ll 1$, so that the additive error term $\Delta$ becomes 
negligibly small compared to the Gaussian term $(2 \pi)^{-d/2} \left( \det B B^T \right)^{-1/2}$.
 For that, we would like the set $Y_i$ to consist of short vectors and to look 
reasonably round. Let us consider the ball $B_r=\left\{x \in {\Bbb R}^n: \ \|x\| \leq r \right\}$ of radius $r$ and choose $Y_i=B_r \cap \Lambda_i$. If the lattice points 
$Y_i$ are sufficiently regular in $B_r \cap \AA_i$ then the moment of inertia of $Y_i$ is roughly the 
moment of inertia of the section $B_r \cap \AA_i$, from which it follows that the maximum 
eigenvalue of $\psi_i$ is about $r^2/\dim \AA_i=r^2/(n-d)$. Roughly,
we get 
$$\rho \approx {r^2 \over (n-d)},$$
where $r$ is the smallest radius of the ball $B_r$ such that the lattice points 
$B_r \cap \Lambda_i$ are distributed regularly in every section $B_r \cap \AA_i$ for 
$i=1, \ldots, d$.

We prove Theorem 2.4 in Section 8.

 In Section 5, we apply Theorem 2.4 to approximate the number of 1-margin multi-way contingency 
 tables, see for example, \cite{Go63} and \cite{DO04}, that is, $\nu$-dimensional 
$k_1 \times \ldots \times k_{\nu}$ arrays of 
non-negative integers $\left(\xi_{j_1 \ldots j_{\nu}} \right)$ with $1 \leq j_i \leq k_i$ for 
$i=1, \ldots, \nu$ with prescribed sums along coordinate hyperplanes $j_i=j$.
 We show that Theorem 2.4 implies that asymptotically the number of such arrays is given by a Gaussian formula 
 (2.3.1)
 as long as $\nu \geq 5$. We suspect that the Gaussian approximation holds as long as 
 $\nu \geq 3$, but the proof would require some additional considerations beyond those of Theorem 2.4.
 
 In \cite{Ba09}, a much cruder asymptotic formula with the main term $e^{g(z)}$ in the logarithmic 
 order is 
 shown to hold for the number of integer points in {\it flow polytopes} (a class of polytopes 
 extending transportation polytopes for $\nu=2$).
 At our request, A. Yong \cite{Yo08} computed a number of examples. Here is one 
 of them, originating in \cite{DE85} and then often used as a  benchmark for various computational approaches:
 \smallskip
 we want to estimate the number of $4 \times 4$ non-negative integer matrices with row sums
 $220, 215, 93$ and 64 and column sums $108, 286, 71$ and 127. The exact number of such matrices
 is $1225914276768514 \approx 1.23 \times 10^{15}$. Framing the problem as the 
 problem of counting integer points in a polytope in the most straightforward way, we obtain an over-determined system $Ax=b$
 (note that the row and column sums of a matrix are not independent). Throwing away one constraint and 
 applying formula (2.3.1), we obtain $1.30 \times 10^{15}$, which overestimates the true number by 
 about $6\%$. The precision is not bad, given that we are applying the Gaussian approximation to the 
 probability mass-function of the sum of 16 independent random 7-dimensional integer vectors.
   
 \subhead (2.5) Gaussian approximation for the number of 0-1 points \endsubhead 
 For a polytope $P$ defined by a system $Ax =b$, $0 \leq x \leq 1$ (shorthand for 
 $0 \leq \xi_j \leq 1$ for $x=\left(\xi_1, \ldots, \xi_n \right)$), we find the point $z=\left(\zeta_1, \ldots, \zeta_n \right)$ maximizing 
 $$h(x)=\sum_{j=1}^n \Bigl( \xi_j \ln {1 \over \xi_j} +(1-\xi_j) \ln {1 \over 1- \xi_j} \Bigr), 
 \qquad x=\left(\xi_1, \ldots, \xi_n \right)$$
 on $P$.  Assuming that $A$ is an integer matrix of rank $d<n$ with the 
 columns $a_1, \ldots, a_n \in {\Bbb Z}^d$,  we compute the $d \times n$ matrix $B$ whose $j$-th column is \break
 $\left(\zeta_j -\zeta_j^2\right)^{1/2} a_j$.
 We approximate the number of 0-1 vectors in $P$ by the Gaussian formula
$$\left| P \cap \{0, 1\}^n \right| \approx {e^{h(z)} \det \Lambda \over (2 \pi)^{d/2} (\det B B^T)^{1/2}},
\tag2.5.1$$
where $\Lambda =A\left({\Bbb Z}^n \right)$. Again, we consider the simplest case of
$\Lambda={\Bbb Z}^d$. We prove the following main result.

\proclaim{(2.6) Theorem} Let us consider a quadratic form $q: {\Bbb R}^d \longrightarrow {\Bbb R}$ 
defined by 
$$q(t) ={1 \over 2} \sum_{j=1}^n \left(\zeta_j -\zeta_j^2 \right) \langle a_j, t \rangle^2.$$
For $i=1, \ldots, d$ let us choose a non-empty finite set $Y_i \subset {\Bbb Z}^n$ such that 
$Ay=e_i$ for all $y \in Y_i$ and let us define a quadratic form 
$\psi_i: {\Bbb R}^n \longrightarrow {\Bbb R}$ by 
$$\psi_i(x)={1 \over |Y_i|} \sum_{y \in Y_i} \langle y, x \rangle^2.$$

Suppose that for some $\lambda >0 $ we have 
$$q(t) \geq \lambda \|t\|^2 \quad \text{for all} \quad t \in {\Bbb R}^d,$$
that for some $\rho >0$ we have 
$$\psi_i(x) \leq \rho \|x\|^2 \quad \text{for all} \quad x \in {\Bbb R}^n \quad \text{and} \quad 
i=1, \ldots, d,$$
that for some $\theta \geq 1$ we have 
$$\|a_j\|_1 \ \leq \ \theta \sqrt{\zeta_j\left(1-\zeta_j\right)} \quad \text{for} \quad j=1, \ldots, n$$
and that for some $0 < \alpha \leq 1/4$
we have 
$$\zeta_j(1-\zeta_j) \geq \alpha \quad \text{for} \quad j=1, \ldots, n.$$
Then, for some absolute constant $\gamma>0$ and for any $0 < \epsilon \leq 1/2$, 
as long as 
$$\lambda \geq \gamma  \epsilon^{-2} \theta^2 \left(d + \ln {1\over \epsilon} \right)^2
\ln \left({n \over \epsilon}\right),$$
we have 
$$\left| P \cap \{0, 1\}^n \right| = e^{h(z)}\left({ \kappa \over (2 \pi)^{d/2} \left( \det B B^T \right)^{1/2}} + 
\Delta \right),$$
where 
$$1-\epsilon \ \leq \kappa \ \leq 1+\epsilon$$
and 
$$|\Delta| \leq \exp\left\{ -{\alpha \over 80 \theta^2 \rho}\right\}.$$ 
\endproclaim

We note that in \cite{Ba08} a much cruder asymptotic formula with the main term $e^{h(z)}$ 
in the logarithmic order is shown 
to hold for the number of 0-1 vectors in flow polytopes.

We prove Theorem 2.6 in Section 7.

In Section 5, we apply Theorem 2.6 to approximate the number of binary 1-margin multi-way contingency 
 tables, see for example, \cite{Go63} and \cite{DO04}, that is, $\nu$-dimensional 
$k_1 \times \ldots \times k_{\nu}$ arrays  $\left(\xi_{j_1 \ldots j_{\nu}} \right)$ 
of $0$'s and $1$'s with $1 \leq j_i \leq k_i$ for 
$i=1, \ldots, \nu$ with prescribed sums along coordinate hyperplanes $j_i=j$.
Alternatively, the number of such arrays is the number of $\nu$-partite uniform hypergraphs 
with prescribed degrees of all vertices. 
 We show that Theorem 2.6 implies that asymptotically the number of such arrays is given by the Gaussian formula 
 (2.5.1)
 as long as $\nu \geq 5$. We suspect that the Gaussian approximation holds as long as 
 $\nu \geq 3$, but the proof would require some additional considerations beyond those of Theorem 2.6.

\head 3. Maximum entropy  \endhead 

We start with the problem of integer point counting.

Let us fix positive numbers $p$ and $q$ such that $p+q=1$.
We recall that a discrete random variable $x$ has geometric distribution if 
$$\Pr \bigl\{x= k \bigr\} =p q^k \quad \text{for} \quad k=0, 1, \ldots.$$
For the expectation and variance of $x$ we have 
$$\EE x ={q \over p} \quad \text{and} \quad \var x={q \over p^2}$$ 
respectively.
Conversely, if $\EE x =\zeta$ for some $\zeta>0$ then 
$$p={1 \over 1+\zeta}, \quad q={\zeta \over 1+\zeta} \quad \text{and} \quad \var x=\zeta+\zeta^2.$$

Our first main result is as follows.
\proclaim{(3.1) Theorem} 
Let $P \subset {\Bbb R}^n$ be the intersection of an affine subspace in ${\Bbb R}^n$ and
the non-negative orthant ${\Bbb R}^n_+$.
Suppose that $P$ is bounded and has a non-empty interior, that is contains a 
point $y=\left(\eta_1, \ldots, \eta_n\right)$ where $\eta_j >0$ for $j=1, \ldots, n$.

Then the strictly concave function 
$$g(x)=\sum_{j=1}^n \Bigl(\left(\xi_j+1 \right) \ln \left(\xi_j+1 \right) -\xi_j \ln \xi_j \Bigr)
\quad \text{for} \quad x=\left(\xi_1, \ldots, \xi_n \right)$$
attains its maximum value on $P$ at a unique point $z=\left(\zeta_1, \ldots, \zeta_n \right)$ such that 
$\zeta_j >0$ for $j=1, \ldots, n$. 

Suppose now that $x_j$ are independent 
geometric random variables with expectations $\zeta_j$ for $j=1, \ldots, n$. 
Let $X=\left(x_1, \ldots, x_n \right)$. Then the probability mass function of $X$ is constant on 
$P \cap {\Bbb Z}^n$ and equal to $e^{-g(z)}$
at every $x \in P \cap {\Bbb Z}^n$. In particular,
$$\left|P \cap {\Bbb Z}^n\right| =e^{g(z)}\Pr\bigl\{ X \in P \bigr\}.$$
\endproclaim 
\demo{Proof} It is straightforward to check that $g$ is strictly concave on the non-negative orthant 
${\Bbb R}^n_+$, so it attains its maximum on $P$ at a unique point 
$z=\left(\zeta_1, \ldots, \zeta_n \right)$. Let us show that $\zeta_j >0$. Since $P$ has a non-empty interior, there is a 
point $y=\left(\eta_1, \ldots, \eta_n \right)$ with $\eta_j >0$ for $j=1, \ldots, n$. We note 
that 
$${\partial \over \partial \xi_j} g =\ln \left({\xi_j +1 \over \xi_j} \right),$$
which is finite for $\xi_j>0$ and equals $+\infty$ for $\xi_j=0$ (we consider the right derivative 
in this case). Therefore, if $\zeta_j=0$ for some $j$ then 
$g\bigl((1-\epsilon) z +\epsilon y \bigr) > g(z)$ for all sufficiently small $\epsilon >0$, which is 
a contradiction. 

Suppose that the affine hull of $P$ is defined by a system of linear equations
$$\sum_{j=1}^n \alpha_{ij} \xi_j =\beta_i \quad \text{for} \quad i=1, \ldots, d.$$
Since $z$ is an interior maximum point, the gradient of $g$ at $z$ is orthogonal to the affine 
hull of $P$, so we have 
$$\ln \left({1 +\zeta_j \over \zeta_j}\right)=\sum_{i=1}^d \lambda_i \alpha_{ij} \quad \text{for}
\quad j=1, \ldots, n$$
and some $\lambda_1, \ldots, \lambda_d$. Therefore, for any $x \in P$, 
$x =\left(\xi_1, \ldots, \xi_n \right)$, we have
$$\sum_{j=1}^n \xi_j \ln \left({1 +\zeta_j \over \zeta_j} \right) =
\sum_{j=1}^n \sum_{i=1}^d \lambda_i \xi_j \alpha_{ij} =\sum_{i=1}^d \lambda_i \beta_i,$$
or, equivalently,
$$\prod_{j=1}^n \left({1 +\zeta_j \over \zeta_j} \right)^{\xi_j} =\exp\left\{ \sum_{i=1}^d \lambda_i \beta_i \right\}. \tag3.1.1$$
Substituting $\xi_j=\zeta_j$ for $j=1, \ldots, n$, we obtain
$$\prod_{j=1}^n \left({1 +\zeta_j \over \zeta_j} \right)^{\zeta_j} =\exp\left\{ \sum_{i=1}^d \lambda_i \beta_i \right\}. \tag3.1.2$$
From (3.1.1) and (3.1.2), we deduce 
$$\split
 \left( \prod_{j=1}^n \left( \zeta_j \over 1+\zeta_j\right)^{\xi_j}\right) \left(\prod_{j=1}^n {1 \over 1+\zeta_j} \right) =&\exp\left\{-\sum_{i=1}^d \lambda_i \beta_i \right\} \left(\prod_{j=1}^n {1 \over 1+\zeta_j} \right)
 \\ =& \prod_{j=1}^n {\zeta_j^{\zeta_j} \over \left(1+\zeta_j\right)^{1+\zeta_j}}= 
 e^{-g(z)}. \endsplit$$
 The last identity states that the probability mass function of $X$ is equal to $e^{-g(z)}$ for every 
 integer point $x \in P$.
 {\hfill \hfill \hfill} \qed
\enddemo

One can observe that the random variable $X$ of Theorem 3.1 has the maximum entropy 
distribution among all distributions on ${\Bbb Z}^n_+$ subject to the constraint $\EE X \in P$.

Theorem 3.1 admits the following straightforward extension.
Let $\ell: {\Bbb R}^n \longrightarrow {\Bbb R}$ be a linear function, 
$$\ell(x)=\gamma_1 \xi_1 + \ldots + \gamma_n \xi_n \quad \text{where} \quad 
x=\left(\xi_1, \ldots, \xi_n \right).$$
Let $P \subset {\Bbb R}^n$ be a polyhedron as in Theorem 3.1, although not necessarily bounded,
and suppose that $\ell$ is bounded on $P$ from above and attains its maximum on $P$ on 
a bounded face of $P$ (it is not hard to see that this condition is sufficient for 
the series $\sum_{x \in P \cap {\Bbb Z}^n} \exp\{\ell(x)\}$ to converge).
 Then the strictly concave function 
$$g_{\ell}(x)=\sum_{j=1}^n \Bigl( \left(\xi_j +1 \right) \ln \left(\xi_j +1 \right) -\xi_j \ln \xi_j +\gamma_j \xi_j \Bigr)$$
attains its maximum on $P$ at a unique point $z=\left(\zeta_1, \ldots, \zeta_n \right)$, where 
$\zeta_j >0$ for $j=1, \ldots, n$. Suppose now that $X=\left(x_1, \ldots, x_n \right)$ is the vector 
of independent geometric random variables such that $\EE x_j =\zeta_j$ for $j=1, \ldots, n$.
Then the probability mass function of $X$ at a point $x \in P \cap {\Bbb Z}^n$ is 
equal to $\exp\left\{-g_{\ell}(z) + \ell(x) \right\}$. In particular, 
$$\sum_{x \in P \cap {\Bbb Z}^n} \exp\{\ell(x)\} \ =\ \exp\left\{g_{\ell}(z) \right\} \Pr\bigl\{X \in P\bigr\}.$$
The proof is a straightforward modification of that of Theorem 3.1.

\subhead (3.2) The Gaussian heuristic for the number of integer points \endsubhead 
Below we provide an informal justification for the Gaussian approximation formula (2.3.1).

Let $P$ be a polytope and let $X$ be 
a random vector as in Theorem 3.1. Suppose that $P$ is defined by a system $Ax=b, x \geq 0$,
where $A=\left(\alpha_{ij} \right)$ is a $d \times n$ matrix of rank $d <n$. Let $Y=AX$, so 
$Y=\left(y_1, \ldots, y_d \right)$, where 
$$y_i=\sum_{j=1}^n \alpha_{ij} x_j \quad \text{for} \quad i=1, \ldots, d.$$
By Theorem 3.1, 
$$\left| P \cap {\Bbb Z}^n\right| =e^{g(z)} \Pr \bigl\{Y = b \bigr\}$$
and
$$\EE Y=A z=b.$$
Moreover, the covariance matrix $Q=\left(q_{ij}\right)$ of $Y$ is computed as follows:
$$q_{ij}=\cov(y_i, y_j) =\sum_{k=1}^n \alpha_{ik} \alpha_{jk} \var x_k= \sum_{k=1}^n \alpha_{ik} \alpha_{jk}
\left( \zeta_k +\zeta_k^2 \right).$$
We would like to approximate the discrete random variable $Y$ by the Gaussian random variable 
$Y^{\ast}$ with the same expectation $b$ and covariance matrix $Q$. We assume now that $A$ is an integer matrix and let $\Lambda =\bigl\{ Ax: \ x \in {\Bbb Z}^n \bigr\}$.
Hence $\Lambda \subset {\Bbb Z}^d$ is a $d$-dimensional lattice. Let $\Pi \subset {\Bbb R}^d$ 
be a fundamental domain of $\Lambda$, so $\vl \Pi=\det \Lambda$. For example, we can choose 
$\Pi$ to be the set of points in ${\Bbb R}^d$ that are closer to the origin than to any other point 
in $\Lambda$. 
Then we can write
$$\left| P \cap {\Bbb Z}^n\right| =e^{g(z)} \Pr \bigl\{Y \in b +\Pi \bigr\}.$$
Assuming that the probability density of $Y^{\ast}$ does not vary much on $b+\Pi$ and that the 
probability mass function of $Y$ at $Y=b$ is well approximated by the integral of the density 
of $Y^{\ast}$ over $b+\Pi$, we obtain (2.3.1).
\bigskip

Next, we consider the problem of counting 0-1 vectors.

Let $p$ and $q$ be positive numbers such that $p+q=1$. We recall that a discrete random variable 
$x$ has Bernoulli distribution if 
$$\Pr\{x=0\}=p \quad \text{and} \quad \Pr\{x=1\}=q.$$
We have 
$$\EE x=q \quad \text{and} \quad \var x= qp.$$
Conversely, if $\EE x=\zeta$ for some $0 < \zeta <1$ then 
$$p=1-\zeta, \quad q=\zeta \quad \text{and} \quad \var x=\zeta -\zeta^2.$$
Our second main result is as follows. 
\proclaim{(3.3) Theorem} Let $P \subset {\Bbb R}^n$ be the intersection of an affine subspace in 
${\Bbb R}^n$ and the unit cube $\bigl\{ 0 \leq \xi_j \leq 1: \quad j=1, \ldots, n \bigr\}$. Suppose that 
$P$ has a non-empty interior, that is, contains a point $y=\left(\eta_1, \ldots, \eta_n \right)$ 
where $0 < \eta_j < 1$ for $j=1, \ldots, n$. Then the strictly concave function 
$$h(x)=\sum_{j=1}^n \left( \xi_j \ln {1 \over \xi_j} + \left(1-\xi_j \right) \ln {1 \over 1-\xi_j} \right)
\quad \text{for} \quad x=\left(\xi_1, \ldots, \xi_n \right)$$
attains its maximum value on $P$ at a unique point $z=\left(\zeta_1, \ldots, \zeta_n \right)$ 
such that $0 < \zeta_j < 1$ for $j=1, \ldots, n$.   

Suppose now that $x_j$ are independent Bernoulli random variables with expectations $\zeta_j$ 
for $j=1, \ldots, n$. Let $X=\left(x_1, \ldots, x_n \right)$. Then the 
probability mass function of $X$ is constant on $P \cap \{0, 1\}^n$ and equal to $e^{-h(z)}$
for every $x \in P \cap \{0, 1\}^n$. In particular,
$$\left| P \cap \{0, 1\}^n \right| \ = \ e^{h(z)}\Pr \bigl\{X \in P \bigr\} .$$
\endproclaim
{\hfill \hfill \hfill} \qed

One can observe that $X$ has the maximum entropy distribution among all distributions on 
$\{0, 1\}^n$ subject to the constraint $\EE X \in P$. 
The proof is very similar to that of Theorem 3.1. Besides, Theorem 3.3 follows from a more general
Theorem 3.5 below.

Again, there is a straightforward extension for exponential sums. For a linear function 
$\ell: {\Bbb R}^n \longrightarrow {\Bbb R}$, 
$$\ell(x)=\gamma_1 \xi_1 + \ldots + \gamma_n \xi_n \quad \text{where} \quad 
x=\left(\xi_1, \ldots, \xi_n \right),$$
we introduce 
$$h_{\ell}(x)=\sum_{j=1}^n \left( \xi_j \ln {1 \over \xi_j} + \left(1-\xi_j \right) \ln {1 \over 1-\xi_j} + \gamma_j
\xi_j\right).$$
Then the maximum value of $h$ on $P$ is attained at a unique point $z=\left(\zeta_1, \ldots, \zeta_n \right)$. If $X=\left(x_1, \ldots, x_n \right)$ is a vector of independent Bernoulli random variables 
such that $\EE x_j =\zeta_j$ then the value of the probability mass function $X$ at a point $ x \in P \cap \{0, 1\}^n$ is equal to $\exp\left\{-h_{\ell}(z) +\ell(x) \right\}$. In particular, 
$$\sum_{x \in P \cap \{0, 1\}^n } \exp\left\{ \ell(x) \right\}\ =\ \exp\left\{h_{\ell}(z) \right\}\Pr\bigl\{ X \in P \bigr\}
 .$$
 
 \subhead (3.4) Comparison with the Monte Carlo method \endsubhead Suppose we want to 
 sample a random 0-1 point from the uniform distribution on $P \cap \{0, 1\}^n$. 
 The standard Monte Carlo rejection method 
 consists in sampling a random 0-1 point $x$, accepting $x$ if $x \in P$ and sampling a new point 
 if $x \notin P$. The probability of hitting $P$ is, therefore, $2^{-n}\left|P \cap {\Bbb Z}^n \right|$.
 It is easy to see that the largest possible value of $h$ in Theorem 3.3 is $n \ln 2$ and is attained 
 at $\zeta_1 =\ldots =\zeta_n=1/2$. Therefore, the rejection sampling using the maximum entropy Bernoulli
 distribution of Theorem 3.3 is at least as efficient as the standard Monte Carlo approach and 
 is essentially more efficient if the value of $h(z)$ is small.

Applying a similar logic as in Section 3.2, we obtain the Gaussian heuristic approximation of 
(2.5.1).

We notice that 
$$h(\xi) =\xi \ln {1 \over \xi} + (1-\xi) \ln {1 \over 1-\xi}$$
is the entropy of the Bernoulli distribution with expectation $\xi$ while 
$$g(\xi)=(\xi+1) \ln (\xi+1) -\xi \ln \xi$$
is the entropy of the geometric distribution with expectation $\xi$. One can suggest the 
following general maximum entropy approach, cf. also a similar computation in \cite{Ja57}.

\proclaim{(3.5) Theorem} Let $S \subset {\Bbb R}^n$ be a finite set and let $\conv(S)$ be 
the convex hull of $S$. Let us assume that $\conv(S)$ has a non-empty interior. 
For $x \in \conv(S)$, let us define $\phi(x)$ to be the maximum 
entropy of a probability distribution on $S$ with expectation $x$, that is, 
$$\split \phi(x)=\max \sum_{s \in S} &p_s \ln {1 \over p_s} \\
\text{Subject to:} \quad \sum_{s \in S} &p_s=1 \\
\sum_{s \in S} s &p_s =x \\
&p_s \geq 0 \quad \text{for all} \quad s \in S. \endsplit$$
Then $\phi(x)$ is a strictly concave continuous function on $\conv(S)$.

Let $A \subset {\Bbb R}^n$ be an affine subspace intersecting the interior of $\conv(S)$.
Then $\phi$ attains its maximum value on $A \cap \conv(S)$ at a unique point $z$
in the interior of $\conv(S)$.
There 
is a unique probability distribution $\mu$ on $S$ with entropy $\phi(z)$ and expectation 
in $A$. Furthermore, the 
probability mass function of $\mu$ is constant on the points of $S \cap A$ and equal to 
$e^{-\phi(z)}$ :
$$\mu\{s\}=e^{-\phi(z)} \quad \text{for all} \quad s \in S \cap A.$$ In particular,
$$\left|S \cap A\right| =e^{\phi(z)}\mu\{S \cap A\}.$$ 
\endproclaim
\demo{Proof} Let 
$$H\Bigl(p_s: \quad s \in S \Bigr)=\sum_{s \in S} p_s \ln {1 \over p_s}$$ be the 
entropy of the probability distribution $\{p_s\}$ on $S$.

Continuity and strict concavity of $\phi$ follows from continuity and strict concavity 
of $H$.  Similarly, uniqueness of $\mu$ follows from the strict concavity of $H$.

Since 
$${\partial \over \partial p_s} H = \ln {1 \over p_s} -1,$$
which is finite for $p_s>0$ and is equal to $+\infty$ for $p_s=0$ (we consider the right derivative),
we conclude that for the optimal distribution $\mu$ we have $p_s>0$ for all $s$. 

Suppose that $A$ is defined by linear equations
$$\langle a_i, x \rangle =\beta_i \quad \text{for} \quad i=1, \ldots, d,$$
where $a_i \in {\Bbb R}^n$ are vectors, $\beta_i \in {\Bbb R}$ are numbers and $\langle \cdot, \cdot \rangle$ is the standard scalar product in ${\Bbb R}^n$. Thus the measure $\mu$ is the solution
to the following optimization problem:
$$\split \sum_{s \in S} &p_s \ln {1 \over p_s} \longrightarrow \max \\ 
\text{Subject to:} \quad \sum_{s \in S} &p_s=1 \\
\sum_{s \in S} \langle a_i, s \rangle &p_s = \beta_i \quad \text{for} \quad i=1, \ldots, d \\
&p_s \geq 0 \quad \text{for all} \quad s \in S. \endsplit$$
Writing the optimality conditions, we conclude that 
for some $\lambda_0, \lambda_1, \ldots, \lambda_d$ we have
$$\ln p_s = \lambda_0 +\sum_{i=1}^d \lambda_i  \langle a_i, s \rangle.$$
Therefore,
$$p_s =  \exp\left\{\lambda_0 + \sum_{i=1}^d \lambda_i \langle a_i, s \rangle \right\}.$$
In particular, for $s \in A$ we have 
$$p_s =\exp\left\{ \lambda_0 + \sum_{i=1}^d \lambda_i \beta_i \right\}.$$
On the other hand,
$$\split \phi(z)=&H\Bigl(p_s: \quad s \in S \Bigr)\\=&-\sum_{s \in S} p_s \left( \lambda_0 +\sum_{i=1}^d \lambda_i \langle a_i, s \rangle\right) \\
=&-\lambda_0 -\sum_{i=1}^d \lambda_i \beta_i,
\endsplit $$
which completes the proof.
{\hfill \hfill \hfill} \qed
\enddemo

Finally, we discuss a continuous version of the maximum entropy approach. 

We recall that $x$ is an exponential random variable with expectation $\zeta>0$ if the density 
function $\psi$ of $x$ is defined by 
$$\psi(\tau)=\cases (1/\zeta) e^{-\tau/\zeta} &\text{for \ } \tau \geq 0 \\ 0 & \text{for\ } \tau <0. \endcases$$
We have 
$$\EE x =\zeta \quad \text{and} \quad \var x=\zeta^2.$$
The characteristic function of $x$ is defined by 
$$\EE e^{i \tau x} ={1 \over 1-i \zeta \tau} \quad \text{for} \quad \tau \in {\Bbb R}.$$
\proclaim{(3.6) Theorem} Let $P \subset {\Bbb R}^n$ be the intersection of an affine subspace in
${\Bbb R}^n$ and a non-negative orthant ${\Bbb R}^n_+$. Suppose that $P$ is bounded and has 
a  non-empty interior. Then the strictly concave function
$$f(x)=n+\sum_{j=1}^n \ln \xi_j \quad \text{for} \quad x=\left(\xi_1, \ldots, \xi_n \right)$$ 
attains its unique maximum on $P$ at a point $z=\left(\zeta_1, \ldots, \zeta_n \right)$, where 
$\zeta_j >0$ for $j=1, \ldots, n$.

Suppose now that $x_j$ are independent exponential random variables with expectations 
$\zeta_j$ for $j=1, \ldots, n$. Let $X=\left(x_1, \ldots, x_n \right)$. Then the  
density of $X$ is constant on $P$ and for every $x \in P$ is equal to $e^{-f(z)}$.
\endproclaim
\demo{Proof} As in the proof of Theorem 3.1, we establish that $\zeta_j>0$ for $j=1, \ldots, n$. 
Consequently, the gradient of $f$ at $z$ must be orthogonal to the affine span of $P$. Assume that 
$P$ is defined by a system of linear equations 
$$\sum_{j=1}^n \alpha_{ij} \xi_j =\beta_i \quad \text{for} \quad i=1, \ldots, d.$$
Then 
$${1 \over \zeta_j} =\sum_{i=1}^d \lambda_i \alpha_{ij} \quad \text{for} \quad j=1, \ldots, n.$$
Therefore, for any $x \in P$, $x=\left(\xi_1, \ldots, \xi_n \right)$, we have 
$$\sum_{j=1}^n {\xi_j \over \zeta_j} =\sum_{i=1}^d \left(\sum_{j=1}^n \alpha_{ij} \xi_j \right)=
\sum_{i=1}^d \lambda_i \beta_i.$$
In particular, substituting $\xi_j =\zeta_j$, we obtain
$$\sum_{j=1}^n {\xi_j \over \zeta_j}=n.$$
Therefore, the density of $X$ at $x \in P$ is equal to 
$$\left( \prod_{j=1}^n {1 \over \zeta_j} \right) \exp\left\{ - \sum_{j=1}^n {\xi_j \over \zeta_j} \right\} =e^{-f(z)}.$$
{\hfill \hfill \hfill} \qed
\enddemo

Again, $X$ has the maximum entropy distribution among all distributions on ${\Bbb R}^n_+$ 
subject to the constraint $\EE X \in P$. 

A similar formula can be obtained for the exponential integral 
$$\int_P e^{\ell(x)} \ dx,$$
where $\ell: {\Bbb R}^n \longrightarrow {\Bbb R}$ is a linear function, 
$$\ell(x)=\gamma_1 \xi_1 + \ldots + \gamma_n \xi_n \quad \text{for} \quad x=\left(\xi_1, \ldots, \xi_n \right).$$
The integral may converge even if $P$ is unbounded. We introduce 
$$f_{\ell}(x) =n+\sum_{j=1}^n \ln \xi_j + \gamma_j \xi_j.$$
If $\ell$ is bounded from above on $P$ and attains its maximum on $P$ on a bounded face then the maximum of $f_{\ell}$ on $P$ is attained at a unique point $z=\left(\zeta_1, \ldots, \zeta_n \right)$.
If $X=\left(x_1, \ldots, x_n \right)$ is a vector of independent exponential random variables such that 
$\EE x_j =\zeta_j$ then the density of $X$ at a point $x \in P$ is equal to 
$$\exp\bigl\{-f_{\ell}(z) + \ell(x) \bigr\}.$$

\subhead (3.7) The Gaussian heuristic for volumes \endsubhead Below we provide an informal 
justification of the Gaussian approximation formula (2.1.1)

Let $P$ be a polytope and let $x_1, \ldots, x_n$ be the random variables as in Theorem 3.6. 
Suppose that $P$ is defined by a system $Ax=b$, $x \geq 0$, where $A=\left(\alpha_{ij}\right)$ is a $d \times n$ 
matrix of rank $d <n$. Let $Y=AX$, so $Y=\left(y_1, \ldots, y_d \right)$, where 
$$y_i =\sum_{j=1}^n \alpha_{ij} x_j \quad \text{for} \quad i=1, \ldots, d.$$  
 In view of Theorem 3.6, the density of $Y$ at $b$ is 
equal to 
$$(\vl P) e^{-f(z)} \left( \det A A^T \right)^{-1/2}$$
(we measure $\vl P$ as the $(n-d)$-dimensional volume with respect to the Euclidean structure 
induced from ${\Bbb R}^n$).

We have $\EE y=b$. The covariance matrix $Q=\left(q_{ij}\right)$ of $Y$ is computed as follows:
$$q_{ij}=\cov \left(y_i, y_j \right)=\sum_{k=1}^n \alpha_{ik} \alpha_{jk} \var x_k =
\sum_{k=1}^n \alpha_{ik} \alpha_{jk} \zeta_k^2.$$
Assuming that the distribution of $Y$ at $Y=b$ is well approximated by the Gaussian distribution, 
we obtain formula (2.1.1)

\head 4. Volumes of multi-index transportation polytopes \endhead

We apply Theorem 2.2 to compute volumes of multi-index transportation polytopes.
We begin our discussion with ordinary (two-index) transportation polytopes. Although Theorem 2.2 
does not imply the validity of the Gaussian approximation here, two-index polytopes provide a simple model case of computations that we later use in the case of 
a larger number of indices.

\subhead (4.1) Transportation polytopes \endsubhead

For integers $m,n >1$ let us choose positive numbers $R=\left(r_1, \ldots, r_m \right)$ and 
$C=\left(c_1, \ldots, c_n\right)$ such that 
$$r_1 + \ldots + r_m = c_1 + \ldots + c_n=N$$ 
and let us consider the polytope $P=P(R,C)$ of all $m \times n$ non-negative 
matrices $x=\left(\xi_{ij} \right)$ with the row sums $r_1, \ldots, r_m$ and the column sums 
$c_1, \ldots, c_n$. As is known, $P$ is a non-empty $(m-1)(n-1)$-dimensional polytope, also 
known as a {\it transportation polytope}, see, for example, \cite{Y+84}.
 If $m=n$ and $R=C=\left(1, \ldots, 1\right)$ then 
$P$ is the polytope of $n \times n$ doubly stochastic matrices, also known as the {\it Birkhoff polytope}.
We note that the row and column sums are not independent, since the total sum of all row sums 
is equal to the total sum of the column sums. We define the affine span of $P$ by the following 
non-redundant system of linear equations:
$$\aligned &\sum_{j=1}^n \xi_{ij} =r_i \quad \text{for} \quad i=1, \ldots, m-1 \\
                    &\sum_{i=1}^m \xi_{ij}=c_j \quad \text{for} \quad j=1, \ldots, n-1 \quad \text{and} \\
                    &\sum \Sb 1 \leq i \leq m \\ 1 \leq j \leq n \endSb \xi_{ij} =N. 
 \endaligned \tag4.1.1$$
In other words, we prescribe the sums of the first $m-1$ rows, the first $n-1$ columns, and the 
total sum of the matrix entries. We observe that every column $a$ of the matrix $A$ of the system 
(4.1.1) contains at most 3 non-zero entries (necessarily equal to 1), so $\|a \| \leq \sqrt{3}$. 

Let $z=\left(\zeta_{ij}\right)$ be a matrix, $z \in P$, maximizing 
$$f(x)=mn + \sum \Sb 1 \leq i \leq m \\ 1 \leq j \leq n \endSb \ln \xi_{ij}.$$
Theorem 2.2 associates with system (4.1.1) the following quadratic form $q$ defined on 
${\Bbb R}^{m+n-1}$:
$$\split q(a; b; \omega) =&{1 \over 2} \sum \Sb 1 \leq i \leq m-1 \\ 1 \leq j \leq n-1 \endSb 
\zeta_{ij}^2 \left(\alpha_i + \beta_j + \omega \right)^2 \\ &\quad + 
{1 \over 2} \sum_{i=1}^{m-1} \zeta_{in}^2 \left(\alpha_i + \omega\right)^2 \\ &\quad + 
{1 \over 2} \sum_{j=1}^{n-1} \zeta_{mj}^2 \left(\beta_j + \omega \right)^2 \\ &\quad + 
{1 \over 2} \zeta_{mn}^2 \omega^2. \endsplit$$
Here 
$$a=\left(\alpha_1, \ldots, \alpha_{m-1} \right), \quad b=\left(\beta_1, \ldots, \beta_{n-1} \right)$$
are real vectors and $\omega$ is a real number, so $(a; b; \omega)$ is interpreted as a vector 
from ${\Bbb R}^{m+n-1}$. 

To bound the eigenvalues of $q$ from below, we bound the eigenvalues of a simpler form 
$$\split \hat{q}(a; b; \omega) = & \sum \Sb 1 \leq i \leq m-1 \\ 1 \leq j \leq n-1 \endSb 
 \left(\alpha_i + \beta_j + \omega \right)^2 \\ &\quad + 
 \sum_{i=1}^{m-1}  \left(\alpha_i + \omega\right)^2 +  \sum_{j=1}^{n-1}  \left(\beta_j + \omega \right)^2  + 
\omega^2. \endsplit$$
Let us consider the $(m-2)$-dimensional subspace $H_a \subset {\Bbb R}^{m+n-1}$ defined by the 
equations
$$\sum_{i=1}^{m-1} \alpha_i =0, \quad b=0, \quad \text{and} \quad \omega =0$$
and the $(n-2)$-dimensional subspace $H_b \subset {\Bbb R}^{m+n-1}$ defined by the 
equations
$$\sum_{i=1}^{n-1} \beta_j =0, \quad a=0, \quad \text{and} \quad \omega =0.$$
We observe that $H_a$ is an eigenspace of $\hat{q}$ with the eigenvalue $n$ 
(since the gradient of $\hat{q}$ at $x \in H_a$ is equal to $2nx$) and that $H_b$ is 
an eigenspace of $\hat{q}$ with the eigenvalue $m$ (since the gradient of $\hat{q}$ at $x \in H_b$ 
is equal to $2mx$). Let $L \subset {\Bbb R}^{m+n-1}$ be the orthogonal complement to 
$H_a +H_b$. Then $\dim L=3$ and $L$ consists of the vectors 
$$\left(\underbrace{\alpha, \ldots, \alpha}_{\text{$m-1$ times}}; 
\underbrace{\beta, \ldots, \beta}_{\text{$n-1$ times}}; \omega \right)$$
for some real $\alpha, \beta$, and $\omega$. Therefore, the restriction of $\hat{q}$ onto $L$ 
can be written as 
$$\split & \hat{q}\left(\underbrace{\alpha, \ldots, \alpha}_{\text{$m-1$ times}}; 
\underbrace{\beta, \ldots, \beta}_{\text{$n-1$ times}}; \omega\right) \\ = 
& \qquad (m-1)(n-1) (\alpha + \beta +\omega)^2+
(m-1)(\alpha+\omega)^2 + (n-1)(\beta + \omega)^2 +\omega^2. \endsplit$$
Since 
$$(\alpha + \beta +\omega)^2 + (\alpha + \omega)^2 +(\beta +\omega)^2 \ \geq \ \delta 
\left( \alpha^2 + \beta^2 +\omega^2 \right)$$
for some absolute constant $\delta>0$ and all $\alpha, \beta$ and $\omega$, we conclude that the eigenvalues of $\hat{q}$ exceed
$$\delta { \min \{m-1, \quad n-1 \}  \over \max\{m-1,\quad  n-1\}}$$
for some absolute constant $\delta >0$. Same holds for the eigenvalues of $q$ as long as 
the numbers $\zeta_{ij}$ are uniformly bounded away from 0. 

We notice that the minimum eigenvalue of $q$ is too small to satisfy the conditions of Theorem 2.2.
In fact, as Canfield and McKay have shown \cite{CM07b}, the volume of the Birkhoff polytope is {\it not} asymptotically 
Gaussian as $m=n \longrightarrow +\infty$, since there is a fourth-order correction akin to the Edgeworth correction. However, a very similar analysis can be applied to certain higher-dimensional versions 
of transportation polytopes and there it produces more satisfying results: asymptotically, volumes of such polytopes turn out to be given by the Gaussian formula (2.1.1). 

\subhead (4.2) Multi-index transportation polytopes \endsubhead
Let us fix an integer $\nu \geq 2$ and let us choose integers $k_1, \ldots, k_{\nu}>1$. 
We consider the 
polytope of $P$ of $k_1 \times \ldots \times k_{\nu}$ arrays of non-negative numbers 
$\xi_{j_1 \ldots j_{\nu}}$, where $1 \leq j_i \leq k_i$ for $i=1, \ldots, \nu$, with prescribed sums along the coordinate hyperplanes. Namely, we choose positive numbers $\beta_{ij}$, where 
 $1 \leq j \leq k_i$ for $i=1, \ldots, \nu$  and such that 
 $$\sum_j \beta_{ij} =N$$ 
 for some $N$ and all $i=1, \ldots, \nu$
 and define $P$ by the inequalities 
 $$\xi_{j_1 \ldots j_{\nu}} \geq 0 \quad \text{for all} \quad j_1, \ldots, j_{\nu}$$ and equations 
$$\aligned &\sum \Sb j_1, \ldots, j_{i-1}, j_{i+1}, \ldots, j_{\nu} \endSb \xi_{j_1 \ldots j_{i-1}, j, j_{i+1} \ldots j_{\nu}}
=\beta_{ij} \\ &\qquad \text{for} \quad i=1, \ldots, \nu \quad \text{and} \quad 1 \leq j \leq k_i-1 \quad 
\text{and} \\
&\sum \Sb j_1, \ldots, j_{\nu} \endSb \xi_{j_1 \ldots j_{\nu}} = N. \endaligned \tag4.2.1$$
Let us choose a pair of indices $1 \leq i \leq \nu$ and $1 \leq j \leq k_i-1$. 
We call the first sum in (4.2.1) the {\it $j$-th sectional sum in direction $i$}.
Hence for each direction $i=1, \ldots, \nu$ we prescribe all but the last one sectional 
sum and also prescribe the total sum of the entries of the array. When $\nu=2$ we obtain the transportation polytope discussed in Section 4.1
We observe that every column $a$ of the matrix $A$ of the system (4.2.1) contains at most $\nu+1$ 
non-zero entries (necessarily equal to 1), so $\|a\| \leq \sqrt{\nu+1}$. 

Let $z=\left(\zeta_{j_1 \ldots j_{\nu}} \right)$ be the point maximizing 
$$f(z) =k_1 \cdots k_{\nu} + \sum \Sb j_1, \ldots, j_{\nu} \endSb \ln \xi_{j_1 \ldots j_{\nu}}$$
on $P$.
We describe the quadratic form $q: {\Bbb R}^d \longrightarrow {\Bbb R}$ which 
Theorem 2.2 associates with system (4.2.1). We have 
$d=k_1 + \ldots + k_{\nu} - \nu+1$ and it is convenient to think of ${\Bbb R}^d$ as 
of a particular coordinate subspace of  a bigger space 
$V={\Bbb R}^{k_1} \oplus \ldots \oplus {\Bbb R}^{k_{\nu}} \oplus {\Bbb R}$. Namely, we think of 
$V$ as of the set of vectors $(t, \omega)$, where 
$$t=\left( \tau_{ij} \right) \quad \text{for} \quad 1 \leq j \leq k_i \quad \text{and} \quad i=1, \ldots, \nu$$
and $\tau_{ij}$ and $\omega$ are real numbers. We identify ${\Bbb R}^d$ with the 
coordinate subspace defined by the equations 
$$\tau_{1 k_1} =\tau_{2 k_2} = \ldots = \tau_{\nu k_{\nu}} =0.$$
Next, we define a quadratic form $p: V \longrightarrow {\Bbb R}$ by 
$$p(t, \omega) = {1 \over 2} \sum \Sb j_1, \ldots, j_{\nu} \endSb \zeta_{j_1 \ldots j_{\nu}}^2 
\left( \tau_{1j_1} + \ldots + \tau_{\nu j_{\nu}} + \omega \right)^2.$$
Then the quadratic form $q$ of Theorem 2.2 is the restriction of $p$ onto ${\Bbb R}^d$.

To bound the eigenvalues of $q$ from below, we consider a simpler quadratic form 
$\hat{q}$ which is the restriction of 
$$\hat{p}(t, \omega) =  \sum \Sb j_1, \ldots, j_{\nu} \endSb 
\left( \tau_{1j_1} + \ldots + \tau_{\nu j_{\nu}} + \omega \right)^2$$
onto ${\Bbb R}^d$. 

For $i=1, \ldots, \nu$, let us consider the $(k_i-2)$-dimensional subspace $H_i \subset {\Bbb R}^d$
defined by the equations
$$\sum_{j=1}^{k_i-1} \tau_{ij}=0, \quad \tau_{i'j} =0 \quad \text{for} \quad i'\ne i \quad 
\text{and all} \quad j, \quad \text{and} \quad \omega =0.$$
Then $H_i$ is an eigenspace of $\hat{q}$ with the eigenvalue 
$$\lambda_i =k_1 \cdots k_{i-1} k_{i+1} \cdots k_{\nu},$$ since the gradient of $\hat{q}$ at $x \in H_i$ is 
equal to $2 \lambda_i x$. Let $L \subset {\Bbb R}^d$ be the orthogonal complement to 
$H_1 \oplus \ldots \oplus H_{\nu}$ in ${\Bbb R}^d$. 
Then $\dim L=\nu +1$ and $L$ 
consists of the vectors 
$$\left(\underbrace{\alpha_1, \ldots, \alpha_1}_{\text{$k_1-1$ times}}, 0;
\underbrace{ \alpha_2, \ldots, \alpha_2}_{\text{$k_2-1$ times}}, 0; \ldots, 
\underbrace{\alpha_{\nu}, \ldots, \alpha_{\nu}}_{\text{$k_{\nu}-1$ times}}; 0; \omega \right)$$ 
for some real $\alpha_1, \ldots, \alpha_{\nu}; \omega$. 
Denoting 
$$\split &\mu_0=\left(k_1-1\right) \cdots \left(k_{\nu}-1\right)  \quad \text{and} \\
&\mu_i =\left(k_1-1\right) \cdots  \left(k_{i-1}- 1 \right) \left(k_{i+1}-1 \right) 
\cdots \left(k_{\nu}-1\right), \endsplit$$
We observe that the restriction of $\hat{q}$
onto $L$ satisfies 
$$\split &\hat{q}\left(\underbrace{\alpha_1, \ldots, \alpha_1}_{\text{$k_1-1$ times}}, 0;
\underbrace{ \alpha_2, \ldots, \alpha_2}_{\text{$k_2-1$ times}}, 0; \ldots, 
\underbrace{\alpha_{\nu}, \ldots, \alpha_{\nu}}_{\text{$k_{\nu}-1$ times}}; 0; \omega \right) \\
& \qquad\geq \  \mu_0 \left(\alpha_1 + \ldots + \alpha_{\nu} +\omega\right)^2 
+\sum_{i=1}^{\nu}  \mu_i \left(\alpha_1 + \ldots + \alpha_{i-1} + \alpha_{i+1} + \ldots +\alpha_{\nu}+\omega \right)^2.
\endsplit$$
Since 
$$\split &\left(\alpha_1 + \ldots + \alpha_{\nu} +\omega\right)^2 + \sum_{i=1}^{\nu}  \left(\alpha_1 + \ldots + \alpha_{i-1} + \alpha_{i+1} + \ldots +\alpha_{\nu}+\omega \right)^2 \\
&\qquad  \geq \ \delta\left(\omega^2 +  \sum_{i=1}^{\nu} \alpha_i^2 \right) \endsplit$$
for some $\delta=\delta(\nu)>0$ and all $\alpha_1, \ldots, \alpha_{\nu}$ and $\omega$,
we conclude that the eigenvalues of $\hat{q}$ exceed
$$\delta(\nu) \min_{i=1, \ldots, \nu} \left(k_i-1\right)^{-2} \prod_{j=1}^{\nu}  \left(k_j-1\right),$$
where $\delta(\nu)>0$ is a constant depending on $\nu$ alone.

Suppose now that $\nu$ is fixed and let us consider a sequence of polytopes $P_n$ 
where $k_1, \ldots, k_{\nu}$ grow roughly proportionately with $n$ and where the coordinates 
$\zeta_{j_1 \ldots j_{\nu}}$ remain in the interval between two positive constants. Then the 
minimum eigenvalue of the quadratic form $q$ in Theorem 2.2 grows as 
$\Omega\left(n^{\nu-2} \right)$. In particular, for $\nu \geq 5$ Theorem 2.2 implies
that the Gaussian formula (2.1.1) approximates the volume of $P_n$ with a relative error which 
approaches 0 as $n$ grows.

As an example, let us consider the (dilated) polytope $P_k$ of {\it polystochastic tensors}, that 
is $k \times \ldots \times k$ arrays of non-negative numbers with all sums along coordinate hyperplanes 
equal to $k^{\nu-1}$, cf. \cite{Gr92}.  By symmetry, we must have 
$$\zeta_{j_1 \ldots j_{\nu}} =1.$$ 
Theorem 2.2 implies that for $\nu \geq 5$
$$\vl P_k= \bigl(1+o(1)\bigr) {e^{k^{\nu}} \over (2 \pi)^{(\nu k -\nu+1)/2}}
\quad \text{as} \quad k \longrightarrow +\infty.$$
Interestingly, for $\nu=2$, where our analysis is not applicable, the formula is smaller by a factor of 
$e^{1/3}$ than the true asymptotic value computed in \cite{CM07b}.

\head 5. The number of multi-way contingency tables \endhead 

We apply Theorems 2.4 and 2.6 to compute the number of multi-way contingency tables.
The smallest eigenvalue of the quadratic form $q$ is bounded as in Section 4 and hence 
our main goal is to bound the additive error $\Delta$.
Again, we begin our discussion with ordinary (two-way) contingency tables, where Theorems 2.4 and 2.6 
do not guarantee the validity of the Gaussian approximation, but which provide a simple model case for computations used later in the 
case of multi-way tables.

\subhead (5.1) Contingency tables \endsubhead 
Let us consider the transportation polytope $P=P(R, C)$, see Section 4.1, where the row sums $r_1, \ldots, r_m$
and the column sums $c_1, \ldots, c_n$ are integer. Integer points in $P(R, C)$ are called 
{\it contingency tables} and 0-1 points in $P(R, C)$ are called {\it binary contingency tables}
with margins $R$ and $C$, see \cite{DE85}.

We assume that 
$P$ is defined by system (4.1.1). To estimate the additive error term $\Delta$ 
in Theorems 2.4 and 2.6, we need to construct 
sets of integer vectors of the following three types: 
\bigskip
for $k=1, \ldots, m-1$ we construct a set 
$Y_k^R$ of $m \times n$ integer matrices $y$ such that the $k$-th row sum of $y$ is 1, all other row 
and column sums, with possible exceptions of the $m$-th row sum and $n$-th column sums are 0,
and the total sum of the matrix entries is 0 as well;
\medskip
for $k=1, \ldots, n-1$ we construct a set $Y_k^C$ of $m \times n$ integer matrices $y$ such that 
the $k$-th column of $y$ sum is 1, all other row and column sums, with possible exceptions of 
the $m$-th row sum and the $n$-th column sums are 0, and the total sum of the matrix entries 
is 0 as well;
\medskip
we construct a set $Y_0$ of $m \times n$ integer matrices $y$ such that all the row and column 
sums of $y$ with possible exceptions of the $m$-th row sum and the $n$-th column sum are 0, and the 
total sum of the matrix entries is 1.
\bigskip
To construct $Y_k^R$, let us choose an integer $1 \leq l \leq n$ and let us 
define a matrix $y=\left(\eta_{ij}\right)$ by letting $\eta_{kl}=1$, $\eta_{ml}=-1$ and letting 
all other entries $\eta_{ij}$ equal to 0. The set $Y_k^R$ contains $n$ vectors $y$ with pairwise 
disjoint support and hence the maximum eigenvalue $\rho_k^R$ of the corresponding quadratic 
form
$$\split \psi^R_k(x) =&{1 \over n} \sum_{y \in Y^R_k} \langle y, x \rangle^2 \\
=&{1 \over n} \sum_{l=1}^n \left( \xi_{kl} -\xi_{ml} \right)^2  \quad \text{for} \quad x=\left(\xi_{ij}\right) \endsplit $$
is $2/n$.
Similarly, to construct $Y_k^C$, let us choose an integer $1 \leq l \leq m$ and let us define 
a matrix $y =\left(\eta_{ij}\right)$ by letting $\eta_{lk}=1$, $\eta_{ln}=-1$ and letting all other entries
$\eta_{ij}$ equal to 0. The maximum eigenvalue $\rho_k^C$ of the corresponding quadratic form
$$\psi^C_k(x) ={1 \over m} \sum_{l=1}^m \left(\xi_{lk}-\xi_{ln} \right)^2$$
is $2/m$.

Finally, to construct $Y_0$, let us choose two indices $1 \leq k \leq m-1$ and $1 \leq l \leq n-1$ and 
let us define a matrix $y=\left(\eta_{ij}\right)$ by letting 
$\eta_{kl}=-1$, $\eta_{kn}=1$, $\eta_{ml}=1$ and letting all other entries $\eta_{ij}$ equal to 0.
For the corresponding quadratic form $\psi_0$ we have 
$$\split \psi_0(x)=&{1 \over (m-1)(n-1)} \sum \Sb 1 \leq k \leq m-1 \\ 1 \leq l \leq n-1 \endSb 
\left(\xi_{kn} +\xi_{ml}-\xi_{kl} \right)^2 \\
\leq & {1 \over (m-1)(n-1)} \sum \Sb 1 \leq k \leq m-1 \\ 1 \leq l \leq n-1 \endSb 
3\left(\xi_{kn}^2 +\xi_{ml}^2+\xi_{kl}^2 \right) \\
= & {3 \over m-1} \sum_{k=1}^{m-1} \xi^2_{kn} + {3 \over n-1} \sum_{l=1}^{n-1} \xi_{ml}^2 
+ {3 \over (m-1)(n-1)} \sum \Sb 1 \leq k \leq m-1 \\ 1 \leq l \leq n-1 \endSb \xi_{kl}^2.  \endsplit $$
Hence we can choose 
$$\rho=\max \left\{{3 \over m-1}, \quad {3 \over n-1} \right\}$$
in Theorems 2.4 and 2.6, so the additive term $\Delta$ is exponentially small in $\min\{m, n\}$. This 
bound is pretty weak but it is getting better as we pass to multiway tables.
In fact, as Canfield and McKay have shown \cite{CM07a}, in the simplest case 
of $R=\left(r, \ldots, r\right)$ and $C=\left(c, \ldots, c \right)$, the number of contingency tables 
is {\it not} given by the Gaussian formula, since there is a 4-th order term correction.

\subhead (5.2) Multi-way contingency tables \endsubhead 
Let us consider the $\nu$-index transportation polytope $P$ of Section 4.2. 
We assume that the affine span of $P$ is defined by system (4.2.1), where numbers 
$\beta_{ij}$ are all integer. The integer points in $P$ are called sometimes {\it multi-way} 
contingency tables while 0-1 points are called binary multi-way contingency tables, see \cite{Go63}
and \cite{DL05}. 

To bound the 
additive error term $\Delta$ in Theorems 2.4 and 2.6, 
we construct a set $Y_{ij}$ of $k_1 \times \ldots \times k_{\nu}$ arrays $y$ of integers such that the total sum of entries of $y$
is 0, the $j$-th sectional sum in the $i$-th direction is 1 all other sectional sums are 0, 
where by ``all other'' we mean all but the $k_i$-th sectional sums in every direction $i=1, \ldots, \nu$.
For that, let us choose $\nu-1$ integers $m_1, \ldots, m_{i-1}, m_{i+1}, \ldots, m_{\nu}$, where
$$1 \leq m_1 \leq k_1, \quad  \ldots,  \quad 1 \leq m_{i-1} \leq k_{i-1},  \quad 1 \leq m_{i+1} \leq k_{i+1}, 
\quad \ldots, \quad 1 \leq m_{\nu} \leq k_{\nu}$$
 and define
$y=\left(\eta_{j_1 \ldots j_{\nu}} \right)$ by letting
$$\eta_{m_1 \ldots m_{i-1}, j, m_{i+1} \ldots m_{\nu}}=1, \quad 
\eta_{m_1 \ldots m_{i-1}, k_i, m_{i+1} \ldots m_{\nu}} =-1$$
and letting all other coordinates of $y$ equal to 0. 

Thus the set $Y_{ij}$ contains $k_1  \cdots k_{i-1} k_{i+1} \cdots k_{\nu}$ elements $y$, and the 
corresponding quadratic form $\psi_{ij}$ can be written as 
$$\split \psi_{ij}(x)=&{1 \over \left|Y_{ij}\right| }
 \sum \Sb m_1, \ldots, m_{i-1}, m_{i+1}, \ldots, m_{\nu} \endSb 
 \left(\xi_{m_1 \cdots m_{i-1}, j, m_{i+1} \cdots m_{\nu}} - \xi_{m_1 \cdots m_{i-1}, k_i, m_{i+1} 
 \cdots m_{\nu}} \right)^2 \\
 &\text{for} \quad x =\left(\xi_{j_1 \ldots j_{\nu}} \right), \endsplit$$
from which the maximum eigenvalue $\rho_{ij}$ of $\psi_{ij}$ is 
$2/k_1 \cdots k_{i-1} k_{i+1} \cdots k_{\nu}$.

Next, we construct a set $Y_0$ of arrays $y$ of $k_1 \cdots k_{\nu}$ integers 
$\left(\eta_{j_1 \ldots j_{\nu}} \right)$ such that the total sum of 
entries of $y$ is $1$ while all sectional sums with a possible exception of the $k_i$-th sectional sum 
in every direction $i$ are equal 0. For that, let us choose $\nu$ integers $m_1, \ldots, m_{\nu}$, where
$$ 1 \leq m_1 \leq k_1-1, \quad \ldots, \quad 1 \leq m_{\nu} \leq k_{\nu}-1$$
 and define
$y=\left(\eta_{j_1, \ldots, j_{\nu}} \right)$ by letting 
$$\split &y_{m_1 \ldots m_{\nu}} = 1-\nu \\
               &y_{k_1, m_2 \ldots m_{\nu}} =1 \\
               &y_{m_1, k_2, m_3 \ldots m_{\nu}} =1 \\
                &\ldots \ldots \ldots \ldots  \\
                &y_{m_1 \ldots m_{\nu-1}, k_{\nu}} =1\endsplit$$
and by letting all other coordinates equal to 0. 

The set $Y_0$ contains 
$\left(k_1-1\right) \cdots \left(k_{\nu}-1 \right)$ elements
 and the corresponding quadratic form $\psi_0$ of Theorems 2.4 and 2.6
can be written as 
$$\split \psi_0(x) = &{1 \over |Y_0|} \sum_{y \in Y_0} \langle y, x \rangle^2 \\ =
&{1 \over |Y_0|} \sum \Sb 1 \leq m_1 \leq k_1-1 \\ \ldots \ldots \ldots \ldots  \\
1 \leq m_{\nu} \leq k_{\nu}-1 \endSb 
\left((1-\nu) \xi_{m_1 \ldots m_{\nu}} + \xi_{k_1, m_2 \ldots m_{\nu}} + \ldots + 
\xi_{m_1 \ldots m_{\nu-1}, k_{\nu}} \right)^2 \\
\leq & {(\nu+1) \over |Y_0|} \sum  \Sb 1 \leq m_1 \leq k_1-1 \\ \ldots \ldots \ldots \ldots  \\
1 \leq m_{\nu} \leq k_{\nu}-1 \endSb 
\left( (1-\nu)^2 \xi_{m_1 \ldots m_{\nu}}^2 + \xi_{k_1, m_2 \ldots m_{\nu}}^2 + \ldots + 
\xi_{m_1 \ldots m_{\nu-1}, k_{\nu}}^2 \right).
\endsplit$$
Therefore, the maximum eigenvalue $\rho_0$ of $\psi_0$ does not exceed 
$$(\nu+1)(\nu-1)^2 \max_{i=1, \ldots, \nu} \left\{ {1 \over \left(k_1-1\right) \cdots \left(k_{i-1}-1\right) 
\left(k_{i+1}-1 \right) \cdots \left(k_{\nu}-1 \right)} \right\}, $$
and the same bound can be used for the value of $\rho$ in Theorems 2.4 and 2.6.

Suppose now that $\nu$ is fixed and let us consider a sequence of polytopes $P_n$ where
$k_1, \ldots, k_{\nu}$ grow roughly proportionately with $n$.
Then in Theorems 2.4 and 2.6 we have
$$\rho=O\left( {1 \over n^{\nu-1}} \right).$$
Let us apply Theorem 2.6 for counting multi-way binary contingency tables. We assume, additionally, 
that for the point $z=\left(\zeta_{j_1 \ldots j_{\nu}} \right)$ maximizing 
$$f(x) = \sum \Sb j_1, \ldots, j_{\nu} \endSb {1 \over \xi_{j_1 \ldots j_{\nu}}} 
\ln {1 \over \xi_{j_1 \ldots j_{\nu}}} +{1 \over 1- \xi_{j_1 \ldots j_{\nu}}}
 \ln {1 \over 1- \xi_{j_1 \ldots j_{\nu}}}$$
on the transportation polytope $P_n$ we have 
$$1-\delta \ \geq \ \zeta_{j_1, \ldots, j_{\nu}}\ \geq \ \delta $$
for some constant $1/2> \delta >0$ and all $j_1, \ldots, j_{\nu}$.
Then we can bound the additive term by 
$$| \Delta| \leq \exp\left\{ -\gamma \delta n^{\nu-1} \right\}$$
for some constant $\gamma >0$.
On the other hand, by Hadamard's inequality,
$$\det B B^{T} =n^{O(n)}.$$
Therefore, for $\nu \geq 3$, the additive term $\Delta$ is negligible compared to the 
Gaussian term. From Section 4.2, we conclude that for $\nu \geq 5$ the relative error for the 
number of multi-way binary contingency tables in $P_n$
for the Gaussian approximation formula (2.5.1) approaches 0 as $n$ grows.

Similarly, we apply Theorem 2.4 for counting multi-way contingency tables. Here we assume,
additionally, that for the point $z=\left(\zeta_{j_1 \ldots j_{\nu}} \right)$ maximizing 
$$f(x) = \sum \Sb j_1, \ldots, j_{\nu} \endSb \left( \xi_{j_1 \ldots j_{\nu}} +1 \right) 
\ln \left( \xi_{j_1 \ldots j_{\nu}} +1 \right) - \xi_{j_1 \ldots j_{\nu}} \ln \xi_{j_1 \ldots j_{\nu}}$$
on the transportation polytope $P_n$ the numbers $\zeta_{j_1 \ldots j_{\nu}}$ lie between 
two positive constants. 
  As in the case of binary tables, we conclude that for $\nu \geq 3$, the additive error term $\Delta$ is negligible compared to the 
Gaussian approximation term as $n \longrightarrow +\infty$. Therefore, for $\nu \geq 5$ the relative error for the 
number of multi-way contingency tables in $P_n$
for the Gaussian approximation formula (2.3.1) approaches 0 as $n$ grows.

Computations show that in the case of $k_1 = \ldots =k_{\nu} =k$ for the matrix $A$ of constraints 
in Theorems 2.4 and 2.6 we have 
$$\det AA^T = k^{(\nu^2 -\nu)(k-1)}.$$
Hence we obtain, for example, that the number of non-negative integer $\nu$-way 
$k \times \ldots \times k$ contingency tables with all sectional sums equal to $r=\alpha k^{\nu-1}$ 
is 
$$\Bigl(1+o(1) \Bigr) \left((\alpha+1)^{\alpha+1}  \alpha^{-\alpha} \right)^{k^{\nu}} 
\left(2 \pi \alpha^2 +2 \pi \alpha \right)^{-(k \nu-\nu+1)/2} k^{(\nu-\nu^2)(k-1)/2}$$
provided $\nu \geq 5$, $k \longrightarrow +\infty$ and $\alpha$ stays between two positive 
constants.
Interestingly, for $\nu=2$ (where our analysis is not applicable) the obtained number is off 
by a constant factor from the true asymptotic obtained in \cite{CM07a}.

Similarly, the number of binary $\nu$-way $k \times \ldots \times k$ binary contingency tables 
with all sectional sums equal to $r=\alpha k^{\nu-1}$ is
$$\Bigl(1+o(1)\Bigr) \left(\alpha^{\alpha} (1-\alpha)^{1-\alpha} \right)^{-k^{\nu}}
\left(2 \pi \alpha - 2 \pi \alpha^2 \right)^{-(k \nu -\nu +1)/2} k^{(\nu-\nu^2)(k-1)/2}$$
as long as $\nu \geq 5$, $k \longrightarrow +\infty$ and $\alpha$ remains separated from $0$ 
and $1$. Again, for $\nu=2$ the formula is off by a constant factor from the asymptotic obtained
in \cite{C+08}.

\head 6. Proof of Theorem 2.2 \endhead 

We start with some standard technical results.
\proclaim{(6.1) Lemma} Let $x_1, \ldots, x_n$ be independent exponential random variables 
such that $\EE x_j =\zeta_j$ for $j=1, \ldots, n$, let $a_1, \ldots, a_n \in {\Bbb R}^d$ be vectors 
which span ${\Bbb R}^d$ and let $Y=x_1 a_1 + \ldots + x_n a_n$. Then the density of $Y$ at 
$b \in {\Bbb R}^d_+$ is equal to
$${1 \over (2 \pi)^d} \int_{{\Bbb R}^d} e^{-i \langle b, t \rangle} 
\left( \prod_{j=1}^n {1 \over 1-i \zeta_j \langle a_j, t \rangle} \right) \ dt.$$
\endproclaim 
\demo{Proof} The characteristic function of $Y$ is
$$\EE e^{i \langle Y, t \rangle} =\prod_{j=1}^n {1 \over 1- i \zeta_j \langle a_j, t \rangle}.$$
The proof now follows by the inverse Fourier transform formula.
{\hfill \hfill \hfill} \qed
\enddemo
We need some standard estimates.
\proclaim{(6.2) Lemma} Let $q: {\Bbb R}^d \longrightarrow {\Bbb R}$ be a positive definite
quadratic form and let $\omega >0$ be a number.
\roster
\item Suppose that $\omega \geq 3$. Then 
$$\int_{t:\  q(t) \geq   \omega d } e^{-q(t)} \ dt  \ \leq \ e^{-\omega d/2} \int_{{\Bbb R}^d} e^{-q(t)} \ dt. $$
\item Suppose that for some $\lambda>0$ we have
$$q(t) \ \geq \ \lambda \|t\|^2 \quad \text{for all} \quad t \in {\Bbb R}^d.$$
Let $a \in {\Bbb R}^d$ be a vector. Then 
$$\int_{t: \ |\langle a, t \rangle| > \omega \|a\|} e^{-q(t)} \ dt \ \leq \ 
e^{-\lambda \omega^2}
\int_{{\Bbb R}^d} e^{-q(t)} \ dt.$$
\endroster
\endproclaim
\demo{Proof} We use the Laplace transform method.
For every $1>\alpha >0$ we have 
$$\split \int_{t:\  q(t) \geq  \omega d} e^{-q(t)} \ dt \ \leq \ 
&\int_{t: \ q(t) \geq  \omega d } \exp\left\{\alpha
\bigl(q(t)- \omega d\bigr) -q(t)
\right\} \ dt \\ \leq\  & e^{-\alpha \omega d} 
\int_{{\Bbb R}^d} \exp\bigl\{-(1-\alpha) q(t) \bigr\} \ dt
\\ =&{e^{-\alpha \omega d} \over (1-\alpha)^{d/2}} \int_{{\Bbb R}^d} e^{-q(t)} \ dt.  \endsplit $$
Optimizing on $\alpha$, we choose $\alpha=1-1/2 \omega$ to conclude that
$$\int_{t: \ q(t) \geq \omega d} e^{-q(t)} \ dt \ \leq \
\exp\left\{ -\omega d +{d \over 2}  +{d \over 2} \ln (2 \omega) \right\} 
 \int_{{\Bbb R}^d} e^{-q(t)} \ dt. $$
 Since 
$$\ln (2\omega ) \leq \omega-1 \quad \text{for} \quad \omega \geq 3,$$
Part (1) follows.

Without loss of generality we assume that $a \ne 0$ in Part (2).
Let us consider the Gaussian probability distribution on ${\Bbb R}^d$ with the density 
proportional to $e^{-q}$. Then $z=\langle a, t \rangle$ is a Gaussian random variable 
such that $\EE z = 0$ and $\var z \leq \|a\|^2/2 \lambda$. Part (2) now follows from the inequality
$$\Pr \bigl\{ |y| \geq \tau \bigr\}  \ \leq \ e^{-\tau^2/2}$$
for the standard Gaussian random variable $y$.
{\hfill \hfill \hfill} \qed
\enddemo

\proclaim{(6.3) Lemma} For $\rho \geq 0$ and $k >d$ we have 
$$\int_{t \in {\Bbb R}^d: \ \|t \| \geq \rho} \left(1+  \|t\|^2 \right)^{-k/2} \ dt 
\ \leq \ {2 \pi^{d/2}  \over \Gamma(d/2) (k-d)} \left(1 + \rho^2 \right)^{(d-k)/2}.$$
\endproclaim
\demo{Proof} Let ${\Bbb S}^{d-1} \subset {\Bbb R}^d$ be the unit sphere in ${\Bbb R}^d$. 
We recall the formula for the surface area of ${\Bbb S}^{d-1}$:
$$\left| {\Bbb S}^{d-1} \right| ={2 \pi^{d/2} \over \Gamma(d/2)}.$$
We have
$$\split \int_{t \in {\Bbb R}^d: \ \|t\| \geq \rho} \left(1+  \|t\|^2 \right)^{-k/2} \ dt =
&\left| {\Bbb S}^{d-1} \right| \int_{\rho}^{+\infty} \left(1 +  \tau^2 \right)^{-k/2} \tau^{d-1} \ d \tau \\
\leq & \left| {\Bbb S}^{d-1} \right| \int_{\rho}^{+\infty} \left(1 + \tau^2 \right)^{(d-k-2)/2} \tau \ d \tau,
\endsplit
$$
where we used that 
$$\tau^{d-1} =\tau \tau^{d-2} \leq \tau \left(1+\tau^2 \right)^{(d-2)/2}.$$
The proof now follows.
{\hfill \hfill \hfill} \qed
\enddemo

Now we are ready to prove Theorem 2.2.

\subhead (6.4) Proof of Theorem 2.2 \endsubhead 
Scaling vectors $a_j$ if necessary, without loss of generality we may assume that 
$\theta =1$. 

From Section 3.7 and Lemma 6.1, we have 
$$\vl P = e^{f(z)} \left( \det A A^T \right)^{1/2} {1 \over (2 \pi)^d} \int_{{\Bbb R}^d} 
e^{-i \langle b, t \rangle} \left( \prod_{j=1}^n {1 \over 1 - i \zeta_j \langle a_j, t \rangle} \right) \ dt.$$
Hence our goal is to estimate the integral and, in particular, to compare it with 
$$\int_{{\Bbb R}^d} e^{-q(t)} \ dt = (2 \pi)^{d/2} \left( \det B B^T \right)^{-1/2}.$$

Let us denote 
$$F(t)=e^{-i \langle b, t \rangle}  \left( \prod_{j=1}^n {1 \over 1 - i \zeta_j \langle a_j, t \rangle} \right) 
\quad \text{for} \quad t \in {\Bbb R}^d.$$
Let 
$$\sigma=4d + 10 \ln {1 \over \epsilon}.$$
We estimate the integral separately over the three regions:
\bigskip
the outer region $\|t\| \geq 1/2$ 
\smallskip
the inner region $q(t) \leq \sigma$ 
\smallskip
 the middle region $\|t\| < 1/2$ and $q(t) > \sigma$.
\bigskip
We note that for a sufficiently large constant $\gamma$ we have $q(t) > \sigma$
in the outer region, we have $\|t\| < 1/2$ in the inner region and the three regions form a 
partition of ${\Bbb R}^d$.
\bigskip
\noindent We start with the outer region $\| t\| \geq 1/2$. Our goal is to show that the integral is negligible
there.

We have 
$$|F(t)|
=\left(\prod_{j=1}^n {1 \over 1 + \zeta_j^2 \langle a_j, t \rangle^2} \right)^{1/2}.$$
Let us denote 
$$\xi_j =\zeta_j^2 \langle a_j, t \rangle^2 \quad \text{for} \quad j=1, \ldots, n.$$
The minimum value of the log-concave function 
$$\prod_{j=1}^n \left(1 + \xi_j \right)$$ 
on the polytope 
$$\sum_{j=1}^n \xi_j \ \geq \ 2\lambda \|t\|^2 \quad \text{and} \quad 0\  \leq \  \xi_j \ \leq \  \|t\|^2$$
is attained at an extreme point of the polytope, that is, at a point where all but possibly one 
coordinate $\xi_j$ is either $0$ or $\|t\|^2$.
Therefore,
$$\left( \prod_{j=1}^n {1 \over 1 + \zeta_j^2 \langle a_j, t \rangle^2} \right)^{1/2} 
\leq \left(1 +  \|t\|^2 \right)^{-\lambda +1/2}. $$
Applying Lemma 6.3, we conclude 
that 
$$\int_{t \in {\Bbb R}^d: \ \|t\| \geq 1/2} |F(t)|  \ dt
\ \leq \ {2 \pi^{d/2} \over \Gamma(d/2)  (2\lambda-d-1) } \left({5 \over 4}\right)^{(d-2\lambda+1)/2}.$$
By the Binet-Cauchy formula and the Hadamard bound, 
$$\det BB^T \leq {n \choose d}  \leq n^d.$$
It follows then that for a sufficiently large absolute constant $\gamma$ and 
the value of the integral over the outer 
region does not exceed 
$(\epsilon /10) (2 \pi)^{d/2} \det (B B^T)^{-1/2}$.

Next, we estimate the integral over the middle region with $\|t\| <1/2$ and $q(t) >\sigma$.
Again, our goal is to show that the integral is negligible.

From the estimate 
$$\left|\ln(1+\xi)-\xi +{\xi^2 \over 2}-{\xi^3 \over 3}   \right|\  \leq \  {|\xi|^4 \over 2} \quad \text{for all complex} \quad |\xi| \leq {1 \over 2},$$
we can write 
$$\ln \left(1 - i \zeta_j \langle a_j, t \rangle \right) =-i \zeta_j \langle a_j, t \rangle 
+{1 \over 2} \zeta_j^2 \langle a_j, t \rangle^2 +{i \over 3} \zeta_j^3 \langle a_j, t \rangle^3 + g_j(t) \zeta_j^4 \langle  a_j, t \rangle^4,$$
where 
$$|g_j(t)| \ \leq \ {1 \over  2} \quad \text{for} \quad j=1, \ldots, n.$$
Since 
$$\sum_{j=1}^n \zeta_j a_j =b,$$
we have 
$$\aligned F(t)=&\exp\left\{ -q(t) - i f(t) +g(t) \right\} \\ 
&\text{where} \quad f(t)={1 \over 3} \sum_{j=1}^n \zeta_j^3 \langle a_j, t \rangle^3
\quad \text{and} \\ & |g(t)| \leq {1 \over 2} \sum_{j=1}^n \zeta_j^4 \langle a_j, t \rangle^4. \endaligned \tag6.4.1$$
In particular,
$$|F(t)| \leq e^{-3q(t)/4}  \quad \text{provided} \quad \|t\| \leq 1/2.$$
Therefore, by Part (1) of Lemma 6.2 we have 
$$\split \left| \ \int\limits \Sb   \|t\| \leq 1/2 \\ q(t)> \sigma  \endSb F(t) \ dt \right| \leq &\int_{t: \ q(t) >\sigma} e^{-3q(t)/4} \ dt  \\
\leq &e^{-3d/2} \epsilon^3 \int_{{\Bbb R}^d} e^{-3q(t)/4} \ dt  \\ \leq &\epsilon^3 \int_{{\Bbb R}^d} e^{-q(t)} \ dt. 
\endsplit $$

Finally, we estimate the integral over the inner region where $q(t) < \sigma$ and,
necessarily, $\|t \| < 1/2$. Here our goal is to show that the integral is very close to 
$\displaystyle \int_{{\Bbb R}^d} e^{-q(t)} \ dt$.

From (6.4.1), we obtain
$$\aligned & \left| \int_{t: \ q(t) < \sigma} F(t) \ dt - \int_{t: \ q(t) < \sigma} e^{-q(t)} \ dt  \right| \\   
&\qquad \qquad \leq 
\int_{t: \ q(t) < \sigma} e^{-q(t)} \left| e^{-i f(t) +g(t)}  -1 \right| \ dt.  \endaligned \tag6.4.2$$
If $q(t) < \sigma$ then $\|t\|^2 \leq \sigma/\lambda$ and hence 
$$|g(t)| \  \leq \ {1 \over 2} \sum_{j=1}^n \zeta_j^4 \langle a_j, t \rangle^4 \ \leq \ 
{\sigma \over2 \lambda} \sum_{j=1}^n \zeta_j^2 \langle a_j, t \rangle^2 = {\sigma^2 \over \lambda}.$$

Thus for all sufficiently large $\gamma$, we have 
$| g(t)| \leq \epsilon/10$.

Let 
$$X=\left\{t : \quad q(t) < \sigma \quad \text{and}  \quad  \zeta_j | \langle a_j, t \rangle| \ \leq {\epsilon \over 10 \sigma} \quad 
\text{for} \quad j=1, \ldots, n  \right\}.$$ By Part (2) of Lemma 6.2, for all sufficiently large $\gamma$, 
we have 
$$\int_{{\Bbb R}^d \setminus X} e^{-q(t)} \ dt  \ \leq \ {\epsilon \over 10}
\int_{{\Bbb R}^d} e^{-q(t)} \ dt$$
whereas for $ t \in X$ we have 
$$|f(t)| \ \leq \  {1 \over 3} \sum_{j=1}^n \zeta_j^3 \left| \langle a_j, t \rangle \right|^3 \ 
\leq \ {\epsilon \over 30 \sigma } \sum_{j=1}^n \zeta_j^2 \langle a_j, t \rangle^2 \ \leq {\epsilon \over 15}.
 $$ 
Estimating 
$$\left|e^{-i f(t) +g(t)} -1 \right| \leq {\epsilon \over 3}  \quad \text{for} \quad t \in X
\quad \text{and} \quad \left|e^{-i f(t) +g(t)} -1 \right| \leq 3 \quad \text{for} \quad t \notin X  $$
we deduce from (6.4.2) that 
$$\split \left| \int_{t: \ q(t) < \sigma} F(t) \ dt - \int_{t: \ q(t) < \sigma} e^{-q(t)} \ dt \right|  \ \leq \ 
&3 \int_{{\Bbb R}^d \setminus X} e^{-q(t)} \ dt + {\epsilon \over 3} \int_{X} e^{-q(t)} \ dt \\ \leq 
& {2 \epsilon \over 3}  \int_{{\Bbb R}^d} e^{-q(t)} \ dt. \endsplit$$
Since by Part (1) of Lemma 6.2, we have 
$$\int_{t: \ q(t) > \sigma} e^{-q(t)} \ dt  \ \leq \ e^{-2d} \epsilon^5 \int_{{\Bbb R}^d} e^{-q(t)} \ dt,$$
the proof follows.
{\hfill \hfill \hfill} \qed

\head 7. Proof of Theorem 2.6 \endhead

First, we represent the number of 0-1 points as an integral.

\proclaim{(7.1) Lemma} Let $p_j, q_j$ be positive numbers such that $p_j+q_j=1$ for $j=1, \ldots, n$
and let $\mu$ be the Bernoulli measure on the set $\{0, 1\}^n$ of 0-1 vectors:
$$\mu\{x\}=\prod_{j=1}^n p_j^{1-\xi_j} q_j^{\xi_j} \quad \text{for} \quad x=\left(\xi_1, \ldots, \xi_n \right).$$
Let $P \subset {\Bbb R}^n$ be a polyhedron defined by a vector equation
$$\xi_1 a_1 + \ldots + \xi_n a_n =b$$
for some integer vectors $a_1, \ldots, a_n; b \in {\Bbb Z}^d$ and 
inequalities
$$0 \ \leq \ \xi_1, \ldots, \xi_n \ \leq \  1.$$
Let $\Pi \subset {\Bbb R}^d$ be the parallelepiped consisting of the points 
$t=\left(\tau_1, \ldots, \tau_d \right)$ such that 
$$-\pi  \  \leq \ \tau_k \ \leq \ \pi  \quad \text{for} \quad k=1, \ldots, d.$$
Then, for 
$$\mu(P)=\sum_{x \in P \cap \{0,1\}^n} \mu\{x\}$$
we have 
$$\mu(P)={1 \over (2 \pi)^d} \int_{\Pi} e^{-i \langle t, b \rangle} \prod_{j=1}^n 
\left(p_j + q_j e^{i \langle a_j, t \rangle} \right) \ dt.$$
Here $\langle \cdot, \cdot \rangle$ is the standard scalar product in ${\Bbb R}^d$ and 
$dt$ is the Lebesgue measure on ${\Bbb R}^d$.
\endproclaim
\demo{Proof} The result follows from the expansion 
$$\prod_{j=1}^n 
\left(p_j + q_j e^{i \langle a_j, t \rangle} \right) =\sum \Sb x \in \{0,1\}^n \\ x=\left(\xi_1, \ldots, \xi_n \right)
\endSb
\exp\left\{i \langle \xi_1 a_1 + \ldots +\xi_n a_n, \ t \rangle \right\} 
\prod_{j=1}^n p_j^{1-\xi_j} q_j^{\xi_j}$$
and the identity  
$${1 \over (2 \pi)^d} \int_{\Pi} e^{i \langle u, t \rangle} \ dt =\cases 1 &\text{if\ } u =0 \\
0 &\text{if \ } u \in {\Bbb Z}^d \setminus \{0\}. \endcases$$
{\hfill \hfill \hfill} \qed
\enddemo

The integrand
$$\prod_{j=1}^n \left(p_j +q_j e^{i \langle a_j, t \rangle}\right)$$
is the characteristic function of $Y=AX$ where $X$ is the multivariate Bernoulli random variable and 
$A$ is the matrix with the columns $a_1, \ldots, a_n$.

The following result is crucial for bounding the additive error $\Delta$.

\proclaim{(7.2) Lemma} Let $A$ be a $d \times n$ integer matrix with the columns 
$a_1, \ldots, a_n \in {\Bbb Z}^d$. For $k=1, \ldots, d$ 
let $Y_k \subset {\Bbb Z}^n$ be a non-empty finite set such that $Ay=e_k$ for 
all $y \in Y_k$, where $e_k$ is the $k$-th standard basis vector. 
Let $\psi_k: {\Bbb R}^n \longrightarrow {\Bbb R}$ be a quadratic form,
$$\psi_k(x) ={1 \over |Y_k|} \sum_{y \in Y_k} \langle y, x \rangle^2 \quad \text{for} \quad
x \in {\Bbb R}^n,$$
and let $\rho_k$ be the maximum eigenvalue of of $\psi_k$.

Suppose further that $0 < \zeta_1, \ldots, \zeta_n <1$ are numbers such that 
$$\zeta_j (1-\zeta_j) \geq \alpha \quad \text{for some} \quad 0 < \alpha \leq 1/4.$$
Then for $t=\left(\tau_1, \ldots, \tau_d \right)$ where $-\pi \leq \tau_k \leq \pi$ for $k=1, \ldots, d$ 
we have 
$$\left| \prod_{j=1}^n \left(1-\zeta_j +\zeta_j e^{i \langle a_j, t \rangle} \right) \right| 
\ \leq \ \exp\left\{-{\alpha \tau_k^2 \over 5\rho_k} \right\}.$$
\endproclaim
\demo{Proof} 
Let us denote
$$F(t)=\prod_{j=1}^n \left( 1-\zeta_j +\zeta_j e^{i \langle a_j, t \rangle}\right).$$
Then
$$|F(t)|^2 =\prod_{j=1}^n \left((1-\zeta_j)^2  +2 \zeta_j (1-\zeta_j) \cos \langle a_j, t \rangle +\zeta_j^2 \right).$$
For real numbers $\xi, \eta$, we write 
$$\xi \equiv \eta \mod 2 \pi$$
if $\xi-\eta$ is an integer multiple of $2 \pi$.
Let 
$$-\pi \leq \gamma_j \leq \pi \quad \text{for} \quad j=1, \ldots, n$$
be numbers such that 
$$\langle a_j, t \rangle \equiv \gamma_j \mod 2 \pi \quad \text{for} \quad j=1, \ldots, n. \tag7.2.1$$
Hence we can write
$$|F(t)|^2 =\prod_{j=1}^n \left((1-\zeta_j)^2 +2 \zeta_j (1-\zeta_j)  \cos \gamma_j +\zeta_j^2 \right).$$
Since 
$$\cos \gamma \ \leq \ 1-{\gamma^2 \over 5} \quad \text{for} \quad -\pi \leq  \gamma \leq \pi,$$
we have
$$|F(t)|^2 \leq \prod_{j=1}^n \left(1 -{2 \zeta_j (1-\zeta_j) \over 5} \gamma_j^2 \right) \ \leq \ 
\exp\left\{ - {2  \alpha \over 5} \sum_{j=1}^n \gamma_j^2 \right\}. \tag7.2.2$$
Let 
$$c=\left(\gamma_1, \ldots, \gamma_n \right), \quad c \in {\Bbb R}^n.$$
Then for all $y \in Y_k$ we have 
$$\tau_k= \langle e_k, t \rangle=\langle Ay, t \rangle =\langle y, A^{\ast} t \rangle \equiv \langle y, c \rangle  \mod 2\pi, $$
where $A^{\ast}$ is the transpose matrix of $A$.
Since $|\tau_k| \leq \pi$, we have 
$$| \langle y, c \rangle | \ \geq \ |\tau_k| \quad \text{for all} \quad y \in Y_k.$$
Therefore,
$$\|c\|^2 \ \geq \ {1 \over \rho_k} \psi_k(c) \ =\ {1 \over \rho_k |Y_k|  }  \sum_{y \in Y_k} \langle y, c 
\rangle^2  \ \geq \ 
  \ {\tau_k^2 \over \rho_k}.$$
The proof follows by (7.2.2).
{\hfill \hfill \hfill} \qed
\enddemo

\subhead (7.3) Proof of Theorem 2.6 \endsubhead
By Theorem 3.3 and Lemma 7.1, we write 
$$\left| P \cap \{0, 1\}^n \right| ={e^{h(z)} \over (2 \pi)^d} \int_{\Pi} e^{-i \langle  b,t \rangle} 
\prod_{j=1}^n \left( 1-\zeta_j + \zeta_j e^{i \langle a_j, t \rangle} \right)\ dt, \tag7.3.1$$
where  $\Pi$ is the parallelepiped consisting of the points $t=\left(\tau_1, \ldots, \tau_d \right)$ 
with $-\pi \leq \tau_k \leq \pi$ for $k=1, \ldots, d$.

Let us denote
$$F(t)=e^{-i \langle b, t \rangle} \prod_{j=1}^n \left( 1 -\zeta_j + \zeta_j e^{i \langle a_j, t \rangle} \right).$$
If 
$$\|t \|_{\infty} \leq {1 \over 4 \theta},$$
we have
$$\left|\langle a_j, t \rangle \right| \leq {1 \over 4} \quad \text{for} \quad j =1, \ldots, n.$$
Using the estimate 
$$\left| e^{i \xi} -1 -i \xi +{\xi^2 \over 2} +i {\xi^3 \over 6} \right| \leq {\xi^4 \over 24} \quad \text{for all real} \quad \xi,$$
we can write
$$\split e^{i \langle a_j, t \rangle}=1+i\langle a_j, t \rangle - {\langle a_j, t \rangle^2 \over 2}
-i {\langle a_j, t \rangle^3 \over 6}
+&g_j(t) \langle a_j, t \rangle^4, \\ 
\text{where} \quad &|g_j(t)| \leq {1 \over 24} \quad \text{for} \quad j=1, \ldots, n. \endsplit$$
Therefore,
$$F(t) =e^{- i\langle b, t \rangle}
\prod_{j=1}^n \left(1+i \zeta_j \langle a_j, t \rangle - \zeta_j {\langle a_j, t \rangle^2 \over 2}
-i\zeta_j {\langle a_j, t \rangle^3 \over 6} + \zeta_j g_j(t) \langle a_j, t \rangle^4 \right).$$
Furthermore, using the estimates 
$$\left| \ln (1 +\xi ) - \xi + {\xi^2 \over 2}-{\xi^3 \over 3} \right| \leq {|\xi|^4 \over 2}
  \quad \text{for all complex} \quad |\xi| \leq 1/2$$
and that
$$\sum_{j=1}^n \zeta_j a_j =b_j,$$
we can write 
$$\split &F(t)=
 e^{ -q(t)  + if(t)+g(t)}, \\
\quad \text{where} \quad
& f(t)={1 \over 6} \sum_{j=1}^n (2\zeta_j -1)  \left( \zeta_j- \zeta_j^2 \right)\langle a_j, t \rangle^3 \quad
\text{and} \\ &|g(t)| \leq 2 \sum_{j=1}^n  \langle a_j, t \rangle^4.
\endsplit \tag7.3.2$$
In particular,
$$|g(t)| \ \leq \ {1 \over 4} q(t) \quad \text{provided} \quad \|t\|_{\infty} \ \leq \ {1 \over 4 \theta}.$$

Let 
$$\sigma=4d + 10 \ln {1\over \epsilon}.$$
We split the integral (7.3.1) over three regions.

The outer region: 
$$\|t\|_{\infty} \geq {1 \over 4 \theta}.$$
We let 
$$\Delta={1 \over (2 \pi)^d} \int\limits \Sb t \in \Pi \\ \|t\|_{\infty} \geq 1/ 4\theta \endSb F(t) \ dt,$$
and use Lemma 7.2 to bound $|\Delta|$.

The middle region: 
$$q(t) \geq \sigma \quad \text{and} \quad \|t\|_{\infty} \leq {1 \over 4 \theta}.$$
From (7.3.2) we obtain 
$$|F(t)| \leq e^{-3q(t)/4}$$ 
and as in the proof of Theorem 2.2 (see Section 6.4), we show that the integral over the region 
is asymptotically negligible for all sufficiently large $\gamma$.

The inner region:
$$q(t) <\sigma.$$ 
Here we have 
$$\|t\|_{\infty} \ \leq \ \|t\| \ \leq \ {\sigma \over \sqrt{\lambda}} \leq {1 \over 4 \theta}$$
provided $\gamma$ is sufficiently large.

If $q(t) < \sigma$ then $\|t\|_{\infty} \leq \|t\| \leq \sqrt{\sigma/\lambda}$ and 
$$|g(t)| \ \leq \ 2 \sum_{j=1}^n \langle a_j, t \rangle^4 \ \leq \ 2\theta^2 
{\sigma \over \lambda} \sum_{j=1}^n \left(\zeta_j - \zeta_j^2 \right) \langle a_j, t \rangle^2
\ \leq 4\ { \theta^2 \sigma^2 \over \lambda}.$$
In particular, if constant $\gamma$ is large enough, we have 
$|g(t)| \leq \epsilon/10$.

As in Section 6.4, we define 
$$X=\left\{ t: \quad q(t) < \sigma \quad \text{and} \quad |\langle a_j, t \rangle | \leq {\epsilon \over 10 \sigma} \quad \text{for} \quad 
j=1, \ldots, n \right\}.$$
Hence for $t \in X$ we have
$$|f(t)| \ \leq \ {1 \over 6} \sum_{j=1}^n \left(\zeta_j -\zeta_j^2 \right)  \left| \langle a_j, t \right \rangle|^3 
\ \leq \ {\epsilon  \over 60 \sigma} \sum_{j=1}^n \left(\zeta_j-\zeta_j^2 \right) \langle a_j, t \rangle^2 
\ \leq \ {\epsilon \over 30}.$$
By Part (2) of Lemma 6.2, for all sufficiently large $\gamma$, we have 
$$\int_{{\Bbb R}^d \setminus X} e^{-q(t)} \ dt  \ \leq \ {\epsilon \over 10} \int_{{\Bbb R}^d} e^{-q(t)} \ dt$$
and the proof is finished as in Section 6.4.
{\hfill \hfill \hfill} \qed

\head 8. Proof of Theorem 2.4 \endhead

We begin with an integral representation for the number of integer points.
\proclaim{(8.1) Lemma} Let $p_j, q_j$ be positive numbers such that $p_j+q_j=1$ for $j=1, \ldots, n$
and let $\mu$ be the geometric measure on the set ${\Bbb Z}^n_+$ of non-negative integer vectors:
$$\mu\{x\}=\prod_{j=1}^n p_j q_j^{\xi_j} \quad \text{for} \quad x=\left(\xi_1, \ldots, \xi_n \right).$$
Let $P \subset {\Bbb R}^n$ be a polyhedron defined by a vector equation
$$\xi_1 a_1 + \ldots + \xi_n a_n =b$$
for some integer vectors $a_1, \ldots, a_n; b \in {\Bbb Z}^d$ and 
inequalities
$$\xi_1, \ldots, \xi_n \geq 0.$$
Let $\Pi \subset {\Bbb R}^d$ be the parallelepiped consisting of the points 
$t=\left(\tau_1, \ldots, \tau_d \right)$ such that 
$$-\pi  \  \leq \ \tau_k \ \leq \ \pi  \quad \text{for} \quad k=1, \ldots, d.$$
Then, for 
$$\mu(P)=\sum_{x \in P \cap {\Bbb Z}^n} \mu\{x\}$$
we have 
$$\mu(P)={1 \over (2 \pi)^d} \int_{\Pi} e^{-i \langle t, b \rangle} \prod_{j=1}^n 
{p_j \over 1-q_j e^{i \langle a_j, t \rangle}} \ dt.$$
Here $\langle \cdot, \cdot \rangle$ is the standard scalar product in ${\Bbb R}^d$ and 
$dt$ is the Lebesgue measure in ${\Bbb R}^d$.
\endproclaim
\demo{Proof} As in the proof of Lemma 7.1, the result follows from the multiple geometric expansion 
$$\prod_{j=1}^n {p_j \over 1-q_j e^{i \langle a_j, t \rangle}} =
\sum \Sb x \in {\Bbb Z}^n_+ \\ x=\left(\xi_1, \ldots, \xi_n \right) \endSb 
\exp\bigl\{ i \langle  \xi_1 a_1 + \ldots +  \xi_n a_n, \ t \rangle \bigr\}\prod_{j=1}^n p_j q_j^{\xi_j}.$$
{\hfill \hfill \hfill} \qed
\enddemo

The integrand 
$$\prod_{j=1}^n {p_j \over 1-q_j e^{i \langle a_j, t \rangle}}$$
is, of course, the characteristic function of $Y=AX$, where $X$ is the multivariate geometric 
random variable and $A$ is the matrix with the columns $a_1, \ldots, a_n$.

The following result is an analogue of Lemma 7.2.

\proclaim{(8.2) Lemma} Let $A$ be a $d \times n$ integer matrix with the columns 
$a_1, \ldots, a_n \in {\Bbb Z}^d$. For $k=1, \ldots, d$ let $Y_k \subset {\Bbb Z}^d$ be a non-empty finite 
set such that $Ay =e_k$ for all $y \in Y_k$, where $e_k$ is the $k$-th standard basis vector in 
${\Bbb Z}^d$. Let $\psi_k: {\Bbb R}^n \longrightarrow {\Bbb R}$ be a quadratic form, 
$$\psi_k(x)={1 \over |Y_k|} \sum_{y \in Y_k} \langle y, x \rangle^2 \quad \text{for} \quad x \in {\Bbb R}^n,$$
and let $\rho_k$ be the maximum eigenvalue of $\psi_k$. 
Suppose further that $\zeta_1, \ldots, \zeta_n  >0$ are numbers such that 
$$\zeta_j (1+\zeta_j) \geq \alpha \quad \text{for some} \quad \alpha >0.$$
Then for $t=\left(\tau_1, \ldots, \tau_d \right)$ where $-\pi \leq \tau_k \leq \pi$ 
for $k=1, \ldots, d$, we have 
$$\left| \prod_{j=1}^n {1 \over 1+\zeta_j -\zeta_j e^{i \langle a_j, t \rangle}} \right| 
\ \leq \  \left(1+{2 \over 5} \alpha \pi^2 \right)^{-m_k} \quad \text{where} \quad 
m_k =\left\lfloor {\tau_k^2 \over \rho_k \pi^2} \right\rfloor.$$
\endproclaim
\demo{Proof} 
Let us denote 
$$F(t) =\prod_{j=1}^n {1 \over 1+\zeta_j -\zeta_j e^{i \langle a_j, t \rangle}}.$$
Then
$$|F(t)|^2 =\prod_{j=1}^n {1 \over 1 + 2 \zeta_j \left(1 + \zeta_j \right) (1 -\cos \langle a_j, t \rangle)} .$$
Let 
$$-\pi \ \leq \ \gamma_j \ \leq \ \pi \quad \text{for} \quad j=1, \ldots, n$$ be numbers such that 
$$\gamma_j \equiv \langle a_j, t \rangle \mod 2\pi \quad \text{for} \quad j=1, \ldots, n.$$
Hence we can write 
$$\split |F(t)|^2 =&\prod_{j=1}^n {1 \over 1 + 2 \zeta_j \left(1 + \zeta_j \right) (1 -\cos \gamma_j)}
\\ \leq &\prod_{j=1}^n {1 \over 1+2 \alpha (1-\cos \gamma_j)}. \endsplit $$
Since
$$\cos \gamma \ \leq \ 1-{\gamma^2 \over 5} \quad \text{for} \quad -\pi \leq \gamma \leq \pi,$$
we estimate
$$|F(t)|^2 \ \leq \ \prod_{j=1}^n \left(1 + {2 \over 5}  \alpha  \gamma_j^2 \right)^{-1}. \tag8.2.1$$
Let 
$$c=\left( \gamma_1, \ldots, \gamma_n \right).$$
As in the proof of Lemma 7.2,
we obtain
$$\|c\|^2 \ \geq \ {\tau_k^2 \over \rho_k}. $$ 
Let us denote $\xi_j=\gamma_j^2$ for $j=1, \ldots, n$. 
The minimum of the log-concave function 
$$\sum_{j=1}^n \ln \left(1 +{2 \over 5} \alpha \xi_j \right)$$
on the polytope defined by the inequalities  $0 \leq \xi_j \leq \pi^2$ for $j =1, \ldots, n$ and 
$$\sum_{j=1}^n \xi_j \ \geq \ {\tau_k^2 \over \rho_k}$$ is attained at an extreme point of the polytope, where all but possibly one coordinate $\xi_j$ is either $0$ or $\pi^2$. The number of non-zero coordinates $\xi_j$ is at least $\tau_k^2 / \rho_k \pi^2$ and the proof follows by (8.2.1).
{\hfill \hfill \hfill} \qed
\enddemo

\subhead (8.3) Proof of Theorem 2.4 \endsubhead  By Theorem 3.1 and Lemma 8.1, we have 
$$\left| P \cap {\Bbb Z}^n \right| ={e^{g(z)} \over (2 \pi)^d} \int_{\Pi} e^{-i \langle t, b \rangle}
\prod_{j=1}^n {1 \over 1 + \zeta_j - \zeta_j e^{i \langle a_j, t \rangle}} \ dt, \tag8.3.1$$
where $\Pi$ is the parallelepiped consisting of the points $t=\left(\tau_1, \ldots, \tau_d \right)$ 
with $-\pi \leq \tau_k \leq \pi$ for $k=1, \ldots, d$. 

Let us denote 
$$F(t)=e^{-i \langle t, b \rangle}
\prod_{j=1}^n {1 \over 1 + \zeta_j - \zeta_j e^{i \langle a_j, t \rangle}}.$$
Similarly to the proof of Theorem 2.6 (see Section 7.3), assuming that $\|t\|_{\infty} \leq 1/4\theta$, 
we write
$$\split & F(t)=
e^{-q(t)  -i f(t) +g(t)}, \\
 \text{where} \quad &f(t)={1 \over 6} \sum_{j=1}^n \left(\zeta_j + \zeta_j^2 \right)(2 \zeta_j+1)
 \langle a_j, t \rangle^3  \quad \text{and} \\
 & |g(t)| \leq 2\sum_{j=1}^n \left(1+\zeta_j\right)^4 \langle a_j, t \rangle^4. \endsplit$$
We let 
$$\sigma = 4 d + 10 \ln {1 \over \epsilon}$$ and 
as in the proof of Theorem 2.6 (see Section 7.3), we split the integral (8.3.1) 
over the three regions:
\bigskip 
the outer region: $\|t\|_{\infty} \geq 1/4 \theta$,
\smallskip
the middle region: $q(t) \geq \sigma$ and $\|t\|_{\infty} \leq 1/4\theta$ and 
\smallskip
the inner region: $q(t) < \sigma$.
\bigskip
For the outer region, we let 
$$\Delta={1 \over (2 \pi)^d} \int \limits\Sb t \in \Pi \\ \|t\|_{\infty} \geq 1/4 \theta \endSb F(t) \ dt$$
and use Lemma 8.2 to bound $\Delta$. 

We have 
$$|F(t)| \ \leq \ e^{-3q(t)/4}$$ 
in the middle region and we bound the integral there as in Section 7.3. 

In the inner region, we have $\|t\|_{\infty} \leq \|t\| \leq \sqrt{\sigma/\lambda}$ and 
$$|g(t)| \ \leq \ 2 \sum_{j=1}^n \left(1 + \zeta_j \right)^4 \langle a_j, t \rangle^4 
\ \leq \ 2{\theta^2 \sigma \over \lambda} \sum_{j=1}^n \left(\zeta_j + \zeta_j^2 \right) \langle a_j, t \rangle^2 \ \leq \  4 { \theta^2 \sigma^2 \over \lambda}.$$
We define 
$$X=\left\{t: \quad q(t) < \sigma \quad \text{and} \quad \left(2\zeta_j+1\right) | \langle a_j, t \rangle| \leq 
{\epsilon \over 10 \sigma} \quad \text{for} \quad j=1, \ldots, n \right\}$$
and note that for $t \in X$ we have 
$$|f(t)| \ \leq \ {1 \over 6} \sum_{j=1}^n \left(2 \zeta_j +1 \right)
\left(\zeta_j + \zeta_j^2 \right) \left| \langle a_j, t \rangle \right|^3 \ 
\leq \ {\epsilon \over 60 \sigma } \sum_{j=1}^n \left(\zeta_j + \zeta_j^2 \right) \langle a_j, t \rangle 
\ \leq \  {\epsilon \over 30}.$$
The proof is finished as in Section 7.3.
{\hfill \hfill \hfill} \qed


\Refs
\widestnumber\key{AAAAA}

\ref\key{Ba97}
\by K. Ball
\paper An elementary introduction to modern convex geometry
\inbook  Flavors of Geometry
\pages 1--58
\bookinfo Math. Sci. Res. Inst. Publ., 31
\publ Cambridge Univ. Press
\publaddr Cambridge
\yr 1997
\endref


\ref\key{Ba08}
\by A. Barvinok
\paper  On the number of matrices and a random matrix with prescribed row and column sums and 0-1 entries
\paperinfo preprint arXiv:0806.1480
\yr 2008
\endref


\ref\key{Ba09}
\by A. Barvinok
\paper Asymptotic estimates for the number of contingency tables, integer flows, and volumes of transportation polytopes 
\jour Int. Math. Res. Notices
\vol 2009
\yr 2009
\pages 348--385
\endref

\ref\key{BV97}
\by M. Brion and M. Vergne
\paper Residue formulae, vector partition functions and lattice points in rational polytopes
\jour J. Amer. Math. Soc. 
\vol 10 
\yr 1997
\pages 797--833
\endref 

\ref\key{C+08}
\by E.R. Canfield, C. Greenhill, and B. McKay
\paper Asymptotic enumeration of dense 0-1 matrices with specified line sums
\jour J. Combin. Theory Ser. A 
\vol 115 
\yr 2008
\pages 32--66
\endref

\ref\key{CM07a}
\by E.R. Canfield and B.D. McKay
\paper  Asymptotic enumeration of contingency tables with constant margins
\paperinfo preprint arXiv math.CO/0703600
\yr 2007
\endref

\ref\key{CM07b}
\by E.R. Canfield and B.D. McKay
\paper The asymptotic volume of the Birkhoff polytope
\paperinfo preprint arXiv:0705.2422
\yr 2007
\endref

\ref\key{DL05}
\by J.A. De Loera
\paper The many aspects of counting lattice points in polytopes \break
\jour  Math. Semesterber. 
\vol 52 
\yr 2005
\pages 175--195
\endref

\ref\key{DO04}
\by J.A. De Loera and S. Onn
\paper The complexity of three-way statistical tables
\jour SIAM J. Comput. 
\vol 33  
\yr 2004
\pages 819--836
\endref

\ref\key{DE85}
\by P. Diaconis and B. Efron
\paper Testing for independence in a two-way table: new interpre- 
tations of the chi-square statistic. With discussions and with a reply by the authors
\jour Ann. Statist. 
\vol 13 
\yr 1985
\pages  845Ð913
\endref

\ref\key{Go63}
\by I.J. Good
\paper Maximum entropy for hypothesis formulation, especially for multidimensional contingency tables
\jour Ann. Math. Statist. 
\vol 34 
\yr 1963 
\pages 911--934
\endref

\ref\key{GK94}
\by P. Gritzmann and V. Klee
\paper On the complexity of some basic problems in computational convexity. II. Volume and mixed volumes
\inbook Polytopes: Abstract, Convex and Computational (Scarborough, ON, 1993)
\pages  373--466
\bookinfo NATO Adv. Sci. Inst. Ser. C Math. Phys. Sci., 440
\publ Kluwer Acad. Publ.
\publaddr Dordrecht
\yr 1994
\endref

\ref\key{Gr92}
\by M.B. Gromova
\paper The Birkhoff-von Neumann theorem for polystochastic matrices 
\jour Selecta Math. Soviet.  
\vol 11  
\yr 1992
\pages 145--158
\endref

\ref\key{Ja57}
\by E.T. Jaynes
\paper Information theory and statistical mechanics
\jour Phys. Rev. (2) 
\vol 106 
\yr 1957
\pages 620--630
\endref

\ref\key{NN94}
\by Yu. Nesterov and A. Nemirovskii
\book Interior-Point Polynomial Algorithms in Convex Programming
\bookinfo SIAM Studies in Applied Mathematics, 13
\publ Society for Industrial and Applied Mathematics (SIAM)
\publaddr  Philadelphia, PA
\yr 1994
\endref

\ref\key{Re88}
\by J.A. Renegar
\paper Polynomial-time algorithm, based on Newton's method, for linear programming
\jour Math. Programming (Ser. A)
\vol 40 
\yr 1988 
\pages 59--93
\endref

\ref\key{Ve05}
\by S. Vempala
\paper Geometric random walks: a survey
\inbook Combinatorial and Computational Geometry
\pages 577--616
\bookinfo Math. Sci. Res. Inst. Publ., 52
\publ Cambridge Univ. Press
\publaddr  Cambridge
\yr 2005
\endref

\ref\key{Y+84}
\by V.A. Yemelichev,  M.M. Koval\"ev and M.K. Kravtsov
\book Polytopes, Graphs and Optimisation
\publ Cambridge University Press
\publaddr Cambridge
\yr 1984
\endref

\ref\key{Yo08}
\by A. Yong
\paper Personal communication
\yr 2008
\endref

\endRefs
\enddocument
\end